\documentclass[12pt, a4paper]{article}
\usepackage[utf8]{inputenc}
\usepackage[T1]{fontenc}
\usepackage{amsmath, amssymb, amsthm}
\usepackage[english]{babel}
\usepackage{hyperref}
\usepackage{xcolor}
\usepackage{geometry}
\usepackage{enumitem}
\usepackage[all]{xy}

\newtheorem{theorem}{Theorem}[section]
\newtheorem{lemma}[theorem]{Lemma}
\newtheorem{proposition}[theorem]{Proposition}
\newtheorem{corollary}[theorem]{Corollary}
\newtheorem{definition}[theorem]{Definition}
\newtheorem{example}[theorem]{Example}
\newtheorem{remark}[theorem]{Remark}

\newcommand{\Hom}{\operatorname{Hom}}
\newcommand{\Tor}{\operatorname{Tor}}
\newcommand{\Ext}{\operatorname{Ext}}
\newcommand{\Spec}{\operatorname{Spec}}
\newcommand{\Max}{\operatorname{Max}}
\newcommand{\Jac}{\operatorname{Jac}}
\newcommand{\Id}{\operatorname{Id}}
\newcommand{\Coker}{\operatorname{Coker}}
\newcommand{\Ker}{\operatorname{Ker}}

\begin{document}
\begin{center}
	{\large  \bf A survey on the uniform $S$-version of rings, modules  and their homological theories}

\vspace{0.5cm}
Xiaolei Zhang, Wei Qi\\

\end{center}
	
	\begin{abstract}
		This survey provides a comprehensive overview of the recent advancements in the theory of ``uniformly \(S\)''-algebraic structures in commutative ring theory. Originating from the classical concepts of Noetherian, coherent, von Neumann regular, and semisimple rings, the introduction of a multiplicative subset \(S\) has led to the development of \(S\)-Noetherian, \(S\)-coherent, and other \(S\)-analogues. However, the element \(s \in S\) in the original definitions often depends on the ideal or module under consideration. To overcome this limitation and enable deeper module-theoretic characterizations, the notion of "uniformly \(S\)" (abbreviated as \(u\)-\(S\)) was introduced. This survey systematically presents the definitions, characterizations, and properties of \(u\)-\(S\)-torsion modules, \(u\)-\(S\)-exact sequences, and the subsequent uniform analogues of fundamental module classes:  \(u\)-\(S\)-finitely presented, \(u\)-\(S\)-Noetherian, \(u\)-\(S\)-coherent, \(u\)-\(S\)-flat, \(u\)-\(S\)-projective, \(u\)-\(S\)-injective, and \(u\)-\(S\)-absolutely pure modules. We then explore the associated uniform homological dimensions, including the \(u\)-\(S\)-weak global dimension, the \(u\)-\(S\)-global dimension, and their interplay with polynomial rings and localizations. The survey also covers structural ring classes such as \(u\)-\(S\)-von Neumann regular, \(u\)-\(S\)-semisimple, \(u\)-\(S\)-Artinian, \(u\)-\(S\)-multiplication rings, and rings with \(u\)-\(S\)-Noetherian spectrum.
	\end{abstract}
{\it Key Words:} Uniformly \(S\)-torsion theory;  Uniformly \(S\)-version of classical rings and modules;  Uniformly \(S\)-homological theory.\\
{\it 2020 Mathematics Subject Classification:} 13D30.

\newpage
\tableofcontents

\newpage

\section{Introduction and Background}

	The study of commutative rings is deeply rooted in finiteness conditions. Noetherian rings, characterized by the ascending chain condition on ideals, form the bedrock of classical commutative algebra. A natural and far-reaching generalization of the Noetherian property was introduced by Anderson and Dumitrescu in their seminal 2002 paper \cite{Anderson2002}. They defined a ring \(R\) to be \(S\)-Noetherian for a given multiplicative subset \(S \subseteq R\) if for every ideal \(I\) of \(R\), there exists a finitely generated subideal \(K \subseteq I\) and an element \(s \in S\) such that \(sI \subseteq K\). This definition elegantly relaxes the requirement that every ideal be finitely generated; instead, each ideal is ``\(S\)-finite'', meaning it is contained in a finitely generated submodule up to multiplication by an element of \(S\). This single paper opened a floodgate of research, leading to the development of \(S\)-analogues for a vast array of ring-theoretic concepts, including \(S\)-Noetherian modules \cite{Bilgin2018}, \(S\)-coherent rings \cite{Bennis2018}, \(c\)-\(S\)-coherent rings \cite{Bennis2018}, \(S\)-Artinian rings \cite{Sengele2020}, \(S\)-GCD domains \cite{Anderson2019}, \(S\)-strong Mori domains \cite{Kim2014}, and \(S\)-almost perfect rings \cite{Bazzoni2019}. These concepts have been instrumental in extending classical results to a broader context and have found applications in various ring constructions such as polynomial rings, power series rings, trivial extensions, and amalgamated algebras \cite{Lim2014, Lim2015, Liu2007}.
	
	A common thread in these initial definitions is that the element \(s \in S\) that witnesses the finiteness condition is not uniform; it can vary depending on the ideal or module in question. For an \(S\)-Noetherian ring \(R\), for instance, for each ideal \(I\) there exists an \(s_I \in S\) with \(s_I I \subseteq K_I\) for some finitely generated \(K_I\). While powerful, this non-uniformity presents significant obstacles when attempting to prove module-theoretic characterizations or ``uniform'' versions of classical theorems, such as the Eakin-Nagata theorem or the Cartan-Eilenberg-Bass theorem, where a single element must control the behavior of all submodules or a whole family of modules.
	
	To address this challenge, the concept of ``uniformity'' was introduced by Zhang in \cite{Zhang2022}. The key idea is to fix a single element \(s \in S\) that works uniformly for all objects in a given class. This led to the definition of \(u\)-\(S\)-torsion modules: an \(R\)-module \(T\) is \(u\)-\(S\)-torsion if there exists an \(s \in S\) such that \(sT = 0\). This simple yet powerful notion paved the way for a systematic ``uniform'' generalization of homological algebra. From this foundation, Zhang and his collaborators, including Qi, Kim, Wang, Chen, Zhao, Hamed, and others, have developed a comprehensive theory of uniformly \(S\)-structures. 
	
	The first major application of this uniform perspective was the introduction of \(u\)-\(S\)-flat modules in \cite{Zhang2022}. A module \(F\) is \(u\)-\(S\)-flat if tensoring with \(F\) preserves \(u\)-\(S\)-exact sequences. This was shown to be equivalent to \(\Tor_1^R(M,F)\) being \(u\)-\(S\)-torsion for all modules \(M\), providing a direct parallel to the classical characterization of flatness. This work also led to the definition of \(u\)-\(S\)-von Neumann regular rings \cite{Zhang2022}, where a single element \(s \in S\) uniformly controls the regularity condition \(sa = ra^2\) for all \(a \in R\). It was proven that a ring is \(u\)-\(S\)-von Neumann regular if and only if every module is \(u\)-\(S\)-flat, a perfect \(u\)-\(S\)-analogue of the classical characterization of von Neumann regular rings.
	
	Following this, Qi et al. \cite{Qi2022} introduced and thoroughly investigated \(u\)-\(S\)-Noetherian rings and modules. A ring \(R\) is \(u\)-\(S\)-Noetherian if there exists a fixed \(s \in S\) such that every ideal \(I\) of \(R\) satisfies \(sI \subseteq K\) for some finitely generated ideal \(K \subseteq I\). This uniform condition allowed for the successful generalization of the Eakin-Nagata-Formanek theorem and, crucially, the Cartan-Eilenberg-Bass theorem. The latter, characterizing \(u\)-\(S\)-Noetherian rings by the property that direct sums of injective modules are \(u\)-\(S\)-injective, was a landmark result demonstrating the power of the uniform approach. In this context, \(u\)-\(S\)-injective modules were defined and studied to characterize \(u\)-\(S\)-Noetherian rings.
	
The notion of Artinian rings has also been tackled from a uniform perspective. Zhang and Qi \cite{Zhang2025a} defined \(u\)-\(S\)-Artinian modules and rings, characterizing them by a descending chain condition that is uniform with respect to \(s\). They proved a fundamental structure theorem: a ring \(R\) is \(u\)-\(S\)-Artinian if and only if it is \(u\)-\(S\)-Noetherian, its \(u\)-\(S\)-Jacobson radical \(\Jac_S(R)\) is \(S\)-nilpotent, and \(R/\Jac_S(R)\) is a \(u\)-\(S/\Jac_S(R)\)-semisimple ring. This theorem beautifully parallels the classical characterization of Artinian rings and demonstrates the deep interplay between the various uniform properties.

	The duality was completed with the introduction of \(u\)-\(S\)-projective modules and \(u\)-\(S\)-semisimple rings by Zhang and Qi in \cite{Zhang2023a}. An \(R\)-module \(P\) is \(u\)-\(S\)-projective if \(\Hom_R(P, -)\) preserves \(u\)-\(S\)-exact sequences, which is equivalent to \(\Ext_R^1(P,M)\) being \(u\)-\(S\)-torsion for all \(M\). A ring is \(u\)-\(S\)-semisimple if every free module is \(u\)-\(S\)-semisimple, leading to a cascade of equivalent conditions, including that every module is \(u\)-\(S\)-projective and every module is \(u\)-\(S\)-injective. 
	
	The homological dimensions associated with these concepts were explored in subsequent works. Zhang \cite{Zhang2023b} defined the \(u\)-\(S\)-weak global dimension of a ring \(R\), denoted \(u\)-\(S\)-w.gl.dim\((R)\), as the supremum of \(u\)-\(S\)-flat dimensions of all \(R\)-modules. This dimension was shown to be related to the classical weak global dimension of the localization \(R_S\). This study also introduced the \(u\)-\(S\)-flat dimension of a module and provided a new local characterization of the classical weak global dimension. In a parallel development, Zhang \cite{Zhang2023c} introduced the \(u\)-\(S\)-global dimension, \(u\)-\(S\)-gl.dim\((R)\), defined as the supremum of \(u\)-\(S\)-projective dimensions (or equivalently, \(u\)-\(S\)-injective dimensions) of all modules. A highlight of this work is the proof that for polynomial rings, under mild conditions, \(u\)-\(S\)-gl.dim\((R[x]) = u\)-\(S\)-gl.dim\((R) + 1\), a striking departure from the classical Hilbert syzygy theorem  and a direct \(S\)-analogue for the weak global dimension.
	
		Further applications of the uniform torsion theory have been made in the study of multiplication modules and rings. Qi and Zhang \cite{Qi2025} introduced \(u\)-\(S\)-multiplication modules and rings, where a single \(s \in S\) uniformly controls the condition that every submodule is of the form \(IM\) up to multiplication by \(s\). This provides a uniform version of the previously studied \(S\)-multiplication modules \cite{Anderson2020, Chhiti2023}. They also explored the behavior of these rings under constructions like idealizations.

	The theory has also been extended to encompass properties of the prime spectrum. Guesmi and Hamed \cite{Guesmi2025} introduced the concept of rings with a \(u\)-\(S\)-Noetherian spectrum, where there exists a uniform \(s \in S\) such that every radical ideal is \(S\)-finitely generated in a radical sense. The authors \cite{Zhang2025c} proved a Hilbert basis theorem for such rings and studied their behavior under various constructions, including flat overrings, Serre's conjecture rings, Nagata rings, Anderson rings, and semigroup rings. This work connects the uniform theory with classical results on rings with Noetherian spectrum \cite{Ohm1968, Ribenboim1985}.

	\section{The Uniformly S-Torsion Theories and Exact Sequences}
	
	The development of uniformly \(S\)-theories begins with a refined notion of torsion. In classical ring theory, an \(S\)-torsion module is one where every element is annihilated by some element of the multiplicative set \(S\). However, this element can depend on the element itself. To achieve uniformity, we require a single element of \(S\) that annihilates the entire module. This leads to the following foundational definition.
	
	\begin{definition}\cite{Zhang2022}
		\label{def:ustorsion}
		Let \(R\) be a ring and \(S\) a multiplicative subset of \(R\). An \(R\)-module \(T\) is called a \textbf{\(u\)-\(S\)-torsion module} (with respect to \(s\)) if there exists an element \(s \in S\) such that \(sT = 0\).
	\end{definition}
	
	This concept is more restrictive than the usual \(S\)-torsion and is crucial for the uniform theory. The following example from \cite{Zhang2022} illustrates the difference.
	
	\begin{example}\cite[Example 2.2]{Zhang2022}
		\label{ex:ustorsionnotstorsion}
		Let \(R = \mathbb{Z}\) and \(S = \{p^n \mid n \ge 0\}\) for a prime \(p\). Consider the \(\mathbb{Z}\)-module \(M = \mathbb{Z}_{(p)}/\mathbb{Z}\). Then:
		\begin{enumerate}
			\item \(M\) is \(S\)-torsion but \textbf{not} \(u\)-\(S\)-torsion.
			\item \(M\) has no maximal \(u\)-\(S\)-torsion submodule.
		\end{enumerate}
	\end{example}
	
	Despite the lack of maximal \(u\)-\(S\)-torsion submodules in general, the following proposition shows that when such a maximal element exists, it is unique.
	
	\begin{proposition}\cite[Proposition 2.4]{Zhang2022}
		\label{prop:uniquemaximal}
		Let \(R\) be a ring and \(S\) a multiplicative subset of \(R\). If an \(R\)-module \(M\) has a maximal \(u\)-\(S\)-torsion submodule, then \(M\) has only one maximal \(u\)-\(S\)-torsion submodule.
	\end{proposition}
	
	The classes of \(u\)-\(S\)-torsion modules and \(S\)-torsion-free modules are intimately related, as shown by the following functorial characterization.
	
	\begin{proposition}\cite[Proposition 2.5]{Zhang2022}
		\label{prop:storsionfreehom}
		Let \(R\) be a ring and \(S\) a multiplicative subset of \(R\). Then an \(R\)-module \(F\) is \(S\)-torsion-free if and only if \(\Hom_R(T,F) = 0\) for any \(u\)-\(S\)-torsion module \(T\).
	\end{proposition}
	
	A direct consequence of this is that the Tor functors applied to a \(u\)-\(S\)-torsion module yield another \(u\)-\(S\)-torsion module.
	
	\begin{corollary}\cite[Corollary 2.6]{Zhang2022}
		\label{cor:torustorsion}
		Let \(R\) be a ring, \(S\) a multiplicative subset of \(R\), and \(T\) a \(u\)-\(S\)-torsion module. Then \(\Tor_n^R(M,T)\) is \(u\)-\(S\)-torsion for any \(R\)-module \(M\) and \(n \ge 0\).
	\end{corollary}
	
	With this refined notion of torsion, we can now define exactness in the uniform sense. The key is that the usual containment conditions are required to hold only after multiplying by a fixed element of \(S\).
	
	\begin{definition}\cite[Definition 2.7]{Zhang2022}  
		\label{def:usexact}
		Let \(R\) be a ring and \(S\) a multiplicative subset of \(R\). Let \(M\), \(N\) and \(L\) be \(R\)-modules.
		\begin{enumerate}
			\item An \(R\)-homomorphism \(f: M \to N\) is called a \textbf{\(u\)-\(S\)-monomorphism} (resp., \textbf{\(u\)-\(S\)-epimorphism}) if \(\Ker(f)\) (resp., \(\Coker(f)\)) is a \(u\)-\(S\)-torsion module.
			\item An \(R\)-homomorphism \(f: M \to N\) is called a \textbf{\(u\)-\(S\)-isomorphism} if it is both a \(u\)-\(S\)-monomorphism and a \(u\)-\(S\)-epimorphism.
			\item An \(R\)-sequence \(M \xrightarrow{f} N \xrightarrow{g} L\) is called \textbf{\(u\)-\(S\)-exact} (at \(N\)) if there is an element \(s \in S\) such that \(s\Ker(g) \subseteq Im(f)\) and \(sIm(f) \subseteq \Ker(g)\).
		\end{enumerate}
	\end{definition}
	
	It is a straightforward exercise to see that a homomorphism is a \(u\)-\(S\)-monomorphism if and only if \(0 \to M \xrightarrow{f} N\) is \(u\)-\(S\)-exact. The properties of \(u\)-\(S\)-torsion modules under \(u\)-\(S\)-exact sequences are fundamental.
	
	\begin{proposition}\cite[Proposition 2.8]{Zhang2022}
		\label{prop:ustorsionexact}
		Let \(R\) be a ring, \(S\) a multiplicative subset of \(R\) and \(M\) an \(R\)-module. Then the following assertions hold.
		\begin{enumerate}
			\item Suppose \(M\) is \(u\)-\(S\)-torsion and \(f:L \to M\) is a \(u\)-\(S\)-monomorphism. Then \(L\) is \(u\)-\(S\)-torsion.
			\item Suppose \(M\) is \(u\)-\(S\)-torsion and \(g:M \to N\) is a \(u\)-\(S\)-epimorphism. Then \(N\) is \(u\)-\(S\)-torsion.
			\item Let \(f:M \to N\) be a \(u\)-\(S\)-isomorphism. If one of \(M\) and \(N\) is \(u\)-\(S\)-torsion, so is the other.
			\item Let \(0 \to L \xrightarrow{f} M \xrightarrow{g} N \to 0\) be a \(u\)-\(S\)-exact sequence. Then \(M\) is \(u\)-\(S\)-torsion if and only if \(L\) and \(N\) are \(u\)-\(S\)-torsion.
		\end{enumerate}
	\end{proposition}
	
	A crucial property of \(u\)-\(S\)-isomorphisms is that they admit a quasi-inverse, making them an equivalence relation in a strong sense.
	
	\begin{lemma}\cite[Lemma 2.1]{Zhang2023a} 
		\label{lem:usisoquasiinverse}
		Let \(R\) be a ring and \(S\) a multiplicative subset of \(R\). Suppose there is a \(u\)-\(S\)-isomorphism \(f:M \to N\) for \(R\)-modules \(M\) and \(N\). Then there is a \(u\)-\(S\)-isomorphism \(g:N \to M\) and \(t \in S\) such that \(f \circ g = t\Id_N\) and \(g \circ f = t\Id_M\).
	\end{lemma}
	
	A remark from \cite{Zhang2023a} clarifies that this is a stronger condition than simply having an isomorphism after localization.
	
	\begin{remark}\cite[Remark 2.2]{Zhang2023a} 
		\label{rem:usisovslociso}
		Let \(R\) be a ring, \(S\) a multiplicative subset of \(R\), and \(M\) and \(N\) be \(R\)-modules. The condition "there is an \(R\)-homomorphism \(f:M \to N\) such that \(f_S:M_S \to N_S\) is an isomorphism" does not imply "there is an \(R\)-homomorphism \(g:N \to M\) such that \(g_S:N_S \to M_S\) is an isomorphism". For instance, take \(R = \mathbb{Z}\), \(S = \mathbb{Z}\setminus\{0\}\), and the inclusion \(f:\mathbb{Z} \hookrightarrow \mathbb{Q}\).
	\end{remark}
	
	A powerful tool in homological algebra is the Five Lemma. Its uniform \(S\)-analogue holds and is often used in diagram chases.
	
	\begin{proposition} \cite[Proposition 1.2]{Zhang2023b}
		\label{prop:usfivelemma}
		Let \(R\) be a ring, \(S\) a multiplicative subset of \(R\). Consider the following commutative diagram with \(u\)-\(S\)-exact rows:
		\[
		\xymatrix@R=0.5cm@C=0.7cm{
			A \ar[r]^{g_1} \ar[d]_{f_A} & B \ar[r]^{g_2} \ar[d]_{f_B} & C \ar[r]^{g_3} \ar[d]_{f_C} & D \ar[r]^{g_4} \ar[d]_{f_D} & E \ar[d]_{f_E} \\
			A' \ar[r]^{h_1} & B' \ar[r]^{h_2} & C' \ar[r]^{h_3} & D' \ar[r]^{h_4} & E'
		}
		\]
		\begin{enumerate}
			\item If \(f_B\) and \(f_D\) are \(u\)-\(S\)-monomorphisms and \(f_A\) is a \(u\)-\(S\)-epimorphism, then \(f_C\) is a \(u\)-\(S\)-monomorphism.
			\item If \(f_B\) and \(f_D\) are \(u\)-\(S\)-epimorphisms and \(f_E\) is a \(u\)-\(S\)-monomorphism, then \(f_C\) is a \(u\)-\(S\)-epimorphism.
			\item If \(f_A\) is a \(u\)-\(S\)-epimorphism, \(f_E\) is a \(u\)-\(S\)-monomorphism, and \(f_B\) and \(f_D\) are \(u\)-\(S\)-isomorphisms, then \(f_C\) is a \(u\)-\(S\)-isomorphism.
			\item If \(f_A, f_B, f_D, f_E\) are all \(u\)-\(S\)-isomorphisms, then \(f_C\) is a \(u\)-\(S\)-isomorphism.
		\end{enumerate}
	\end{proposition}
	
	The uniform \(S\)-theory also interacts well with the classical concept of \(S\)-finiteness. Recall that a module \(M\) is \textbf{\(S\)-finite} (with respect to \(s\)) if there exists a finitely generated submodule \(F\) and an element \(s \in S\) such that \(sM \subseteq F\).
	
	\begin{proposition}\cite[Proposition 2.9]{Zhang2022}
		\label{prop:sfiniteness}
		Let \(R\) be a ring, \(S\) a multiplicative subset of \(R\) and \(M\) an \(R\)-module. Then the following assertions hold.
		\begin{enumerate}
			\item Let \(M\) be an \(S\)-finite \(R\)-module and \(f:M \to N\) a \(u\)-\(S\)-epimorphism. Then \(N\) is \(S\)-finite.
			\item Let \(0 \to L \xrightarrow{f} M \xrightarrow{g} N \to 0\) be a \(u\)-\(S\)-exact sequence. If \(L\) and \(N\) are \(S\)-finite, so is \(M\).
			\item Let \(f:M \to N\) be a \(u\)-\(S\)-isomorphism. If one of \(M\) and \(N\) is \(S\)-finite, so is the other.
		\end{enumerate}
	\end{proposition}
	
	These foundational tools—\(u\)-\(S\)-torsion, \(u\)-\(S\)-exactness, and the behavior of \(S\)-finiteness—form the bedrock upon which all subsequent uniform concepts are built. They allow for the development of a robust homological algebra where classical proofs can be adapted by replacing exactness with \(u\)-\(S\)-exactness and zero modules with \(u\)-\(S\)-torsion modules.
	
	\section{Finiteness conditions in uniformly S-torsion theories}
	
	Finiteness conditions are the heart of commutative algebra. The \(S\)-Noetherian property was introduced to relax the classical Noetherian condition. Its uniform counterpart, the \(u\)-\(S\)-Noetherian property, strengthens it by requiring the same element \(s \in S\) to work for all ideals. This uniform version enables powerful module-theoretic characterizations that were elusive for the original \(S\)-Noetherian rings.
	
	\subsection{Uniformly S-Finitely Presented Modules}
	
	Before defining \(u\)-\(S\)-Noetherian rings, we need a uniform version of finite presentation, which synthesizes the concepts of \(S\)-finitely presented and \(c\)-\(S\)-finitely presented modules.
	
	\begin{definition}\cite[Definition 2.1]{Zhang2025b}
		\label{def:usfinitelypresented}
		Let \(R\) be a ring, \(S\) be a multiplicative subset of \(R\) and \(s \in S\). An \(R\)-module \(M\) is called \textbf{\(u\)-\(S\)-finitely presented} (with respect to \(s\)) provided that there is an exact sequence
		\[
		0 \to T_1 \to F \xrightarrow{f} M \to T_2 \to 0
		\]
		with \(F\) finitely presented and \(sT_1 = sT_2 = 0\).
	\end{definition}
	
	This definition means that a module is \(u\)-\(S\)-finitely presented if it is \(u\)-\(S\)-isomorphic to a finitely presented module. The class of such modules enjoys good closure properties.
	
	\begin{proposition}\cite[Proposition 2.2]{Zhang2025b}
		\label{prop:usfpclosure}
		Let \(\Phi : 0 \to M \xrightarrow{f} N \xrightarrow{g} L \to 0\) be a \(u\)-\(S\)-exact sequence of \(R\)-modules. The following statements hold.
		\begin{enumerate}
			\item The class of \(u\)-\(S\)-finitely presented modules is closed under \(u\)-\(S\)-isomorphisms.
			\item If \(M\) and \(L\) are \(u\)-\(S\)-finitely presented, so is \(N\).
			\item Any finite direct sum of \(u\)-\(S\)-finitely presented modules is \(u\)-\(S\)-finitely presented.
			\item If \(N\) is \(u\)-\(S\)-finitely presented, then \(L\) is \(u\)-\(S\)-finitely presented if and only if \(M\) is \(S\)-finite.
		\end{enumerate}
	\end{proposition}
	
	The following proposition connects \(u\)-\(S\)-finitely presented modules with the classical notion of finite presentation over the localization \(R_S\), under the mild condition that \(S\) consists of finitely many elements.
	
	\begin{proposition}\cite[Proposition 2.4]{Zhang2025b}
		\label{prop:usfploc}
		Let \(R\) be a ring, \(S\) a multiplicative subset of \(R\) consisting of finite elements. Then an \(R\)-module \(M\) is a \(u\)-\(S\)-finitely presented \(R\)-module if and only if \(M_S\) is a finitely presented \(R_S\)-module.
	\end{proposition}
	
	A local-global principle for finite generation is a classic result. Its uniform version is a precursor to characterizations of Noetherian modules.
	
	\begin{lemma} \cite[Lemma 2.5]{Zhang2025b}
		\label{lem:finitegenlocal}
		Let \(R\) be a ring, \(S\) be a multiplicative subset of \(R\) and \(M\) be an \(R\)-module. The following statements are equivalent:
		\begin{enumerate}
			\item \(M\) is a finitely generated \(R\)-module;
			\item \(M\) is \(\mathfrak{p}\)-finite for any \(\mathfrak{p} \in \Spec(R)\);
			\item \(M\) is \(\mathfrak{m}\)-finite for any \(\mathfrak{m} \in \Max(R)\).
		\end{enumerate}
	\end{lemma}
	
	A similar local-global principle holds for finite presentation, which provides a new characterization of classically finitely presented modules.
	
	\begin{proposition}\cite[Proposition 2.6]{Zhang2025b} 
		\label{prop:finiteplocal}
		Let \(R\) be a ring and \(M\) be an \(R\)-module. The following statements are equivalent:
		\begin{enumerate}
			\item \(M\) is a finitely presented \(R\)-module;
			\item \(M\) is \(u\)-\(\mathfrak{p}\)-finitely presented for any \(\mathfrak{p} \in \Spec(R)\);
			\item \(M\) is \(u\)-\(\mathfrak{m}\)-finitely presented for any \(\mathfrak{m} \in \Max(R)\).
		\end{enumerate}
	\end{proposition}
	
		\subsection{Uniformly S-Noetherian Rings and Modules}
	
	The notion of \(u\)-\(S\)-finiteness can be extended to families of modules. A family of modules is uniformly \(S\)-finite if the same element \(s \in S\) works for all members. This concept is key to defining \(u\)-\(S\)-Noetherian modules.
	
	\begin{definition}\cite[Definition 2.7]{Qi2022}
		\label{def:usfinitefamily}
		Let \(\{M_j\}_{j \in \Gamma}\) be a family of \(R\)-modules and let \(N_j\) be a submodule of \(M_j\) generated by \(\{m_{i,j}\}_{i \in \Lambda_j} \subseteq M_j\) for each \(j \in \Gamma\). The family is said to be \textbf{\(u\)-\(S\)-finite} (with respect to \(s\)) if there exists \(s \in S\) such that \(sM_j \subseteq N_j\) for each \(j \in \Gamma\).
	\end{definition}
	
	An \(R\)-module \(M\) is called a \(u\)-\(S\)-Noetherian module if the set of all its submodules is \(u\)-\(S\)-finite.
	
	\begin{definition}\cite[Definition 2.7]{Qi2022}
		\label{def:usnoetherianmodule}
		Let \(R\) be a ring and \(S\) a multiplicative subset of \(R\). An \(R\)-module \(M\) is called a \textbf{uniformly \(S\)-Noetherian} (or \(u\)-\(S\)-Noetherian) \(R\)-module provided that the set of all submodules of \(M\) is \(u\)-\(S\)-finite.
	\end{definition}
	
	This definition leads to the following elegant characterization in terms of a uniform ascending chain condition.
	
	\begin{theorem}\cite[Theorem 2.8]{Qi2022}
		\label{thm:usnoetherianacc}
		(Eakin-Nagata-Formanek Theorem for uniformly \(S\)-Noetherian modules)
		Let \(R\) be a ring, \(S\) a multiplicative subset of \(R\), and \(M\) an \(R\)-module. Then the following conditions are equivalent:
		\begin{enumerate}
			\item \(M\) is uniformly \(S\)-Noetherian;
			\item there exists \(s \in S\) such that any ascending chain of submodules of \(M\) is stationary with respect to \(s\);
			\item there exists \(s \in S\) such that any nonempty subset of submodules of \(M\) has a maximal element with respect to \(s\).
		\end{enumerate}
	\end{theorem}
	
	Applying this to the ring itself gives a characterization of \(u\)-\(S\)-Noetherian rings.
	
	\begin{corollary}\cite[Corollary 2.9]{Qi2022}
		\label{cor:usnoetherianringacc}
		Let \(R\) be a ring and \(S\) a multiplicative subset of \(R\). Then the following conditions are equivalent:
		\begin{enumerate}
			\item \(R\) is uniformly \(S\)-Noetherian;
			\item there exists \(s \in S\) such that any ascending chain of ideals of \(R\) is stationary with respect to \(s\);
			\item there exists \(s \in S\) such that any nonempty subset of ideals of \(R\) has a maximal element with respect to \(s\).
		\end{enumerate}
	\end{corollary}
	
	When \(S\) satisfies the maximal multiple condition (i.e., there exists \(s \in S\) such that every element of \(S\) divides \(s\)), the \(S\)-Noetherian and \(u\)-\(S\)-Noetherian properties coincide.
	
	\begin{proposition} \cite[Proposition 2.5]{Qi2022}
		\label{prop:usnoethiffnoeth}
		Let \(R\) be a ring and \(S\) a multiplicative subset of \(R\) satisfying the maximal multiple condition. Then the following statements are equivalent.
		\begin{enumerate}
			\item \(R\) is a uniformly \(S\)-Noetherian ring (resp., uniformly \(S\)-PIR).
			\item \(R\) is an \(S\)-Noetherian ring (resp., \(S\)-PIR).
			\item \(R_S\) is a Noetherian ring (resp., PIR).
		\end{enumerate}
	\end{proposition}

	The \(u\)-\(S\)-Noetherian property behaves well under \(u\)-\(S\)-exact sequences.
	
	\begin{lemma} \cite[Lemma 2.13]{Qi2022}
		\label{lem:usnoetherianexact}
		Let \(R\) be a ring and \(S\) a multiplicative subset of \(R\). Let \(0 \to A \to B \to C \to 0\) be an exact sequence of \(R\)-modules. Then \(B\) is uniformly \(S\)-Noetherian if and only if \(A\) and \(C\) are uniformly \(S\)-Noetherian.
	\end{lemma}
	
	This extends to \(u\)-\(S\)-exact sequences as well.
	
	\begin{proposition} \cite[Proposition 2.14]{Qi2022}
		\label{prop:usnoetherianusexact}
		Let \(R\) be a ring and \(S\) a multiplicative subset of \(R\). Let \(0 \to A \to B \to C \to 0\) be a \(u\)-\(S\)-exact sequence of \(R\)-modules. Then \(B\) is uniformly \(S\)-Noetherian if and only if \(A\) and \(C\) are uniformly \(S\)-Noetherian.
	\end{proposition}
	
	A direct consequence is that the property is invariant under \(u\)-\(S\)-isomorphisms.
	
	\begin{corollary}\cite[Corollary 2.15]{Qi2022}
		\label{cor:usnoetherianisoinv}
		Let \(R\) be a ring, \(S\) a multiplicative subset of \(R\), and \(M \xrightarrow{f} N\) a \(u\)-\(S\)-isomorphism of \(R\)-modules. If one of \(M\) and \(N\) is uniformly \(S\)-Noetherian, then so is the other.
	\end{corollary}
	
	This leads to a new local characterization of classical Noetherian modules.
	
	\begin{proposition}\cite[Proposition 2.16]{Qi2022}
		\label{prop:noetherianlocal}
		Let \(R\) be a ring and \(M\) an \(R\)-module. Then the following conditions are equivalent:
		\begin{enumerate}
			\item \(M\) is Noetherian;
			\item \(M\) is uniformly \(\mathfrak{p}\)-Noetherian for any \(\mathfrak{p} \in \Spec(R)\);
			\item \(M\) is uniformly \(\mathfrak{m}\)-Noetherian for any \(\mathfrak{m} \in \Max(R)\).
		\end{enumerate}
	\end{proposition}
	
	The following theorem gives several equivalent characterizations of \(u\)-\(S\)-Noetherian rings, including a new one in terms of \(u\)-\(S\)-finite presentation of finitely generated modules.
	
	\begin{theorem}\cite[Theorem 2.7]{Zhang2025b}
		\label{thm:usnoetherianchar}
		Let \(R\) be a ring and \(S\) be a multiplicative subset of \(R\). Then the following statements are equivalent:
		\begin{enumerate}
			\item A ring \(R\) is \(u\)-\(S\)-Noetherian;
			\item Any \(S\)-finite module is \(u\)-\(S\)-Noetherian;
			\item Any finitely generated module is \(u\)-\(S\)-Noetherian;
			\item There is \(s \in S\) such that any finitely generated module is \(u\)-\(S\)-finitely presented with respect to \(s\).
		\end{enumerate}
	\end{theorem}

	\subsection{Uniformly S-Coherent Rings and Modules}
	
	Coherent rings are a natural generalization of Noetherian rings where the finiteness condition is imposed on the presentation of finitely generated ideals. The uniform \(S\)-analogue is defined as follows.
	
	\begin{definition} \cite[Definition 3.1]{Zhang2025b}
		\label{def:uscoherentmodule}
		Let \(R\) be a ring and \(S\) be a multiplicative subset of \(R\). An \(R\)-module \(M\) is called a \textbf{\(u\)-\(S\)-coherent module} (with respect to \(s\)) provided that there is \(s \in S\) such that it is \(S\)-finite with respect to \(s\) and any finitely generated submodule of \(M\) is \(u\)-\(S\)-finitely presented with respect to \(s\).
	\end{definition}
	
	The class of \(u\)-\(S\)-coherent modules is well-behaved under \(u\)-\(S\)-exact sequences.
	
	\begin{theorem}\cite[Theorem 3.2]{Zhang2025b}
		\label{thm:uscoherentexact}
		Let \(\Phi : 0 \to M \xrightarrow{f} N \xrightarrow{g} L \to 0\) be a \(u\)-\(S\)-exact sequence of \(R\)-modules. The following statements hold.
		\begin{enumerate}
			\item The class of \(u\)-\(S\)-coherent modules is closed under \(u\)-\(S\)-isomorphisms.
			\item If \(L\) is \(u\)-\(S\)-coherent, then \(M\) is \(u\)-\(S\)-coherent if and only if \(N\) is \(u\)-\(S\)-coherent.
			\item Any finite direct sum of \(u\)-\(S\)-coherent modules is \(u\)-\(S\)-coherent.
			\item If \(N\) is \(u\)-\(S\)-coherent and \(M\) is \(S\)-finite, then \(L\) is \(u\)-\(S\)-coherent.
		\end{enumerate}
	\end{theorem}
	
	A ring is called \(u\)-\(S\)-coherent if it is \(u\)-\(S\)-coherent as a module over itself.
	
	\begin{definition}\cite[Definition 3.5]{Zhang2025b}
		\label{def:uscoherentring}
		Let \(R\) be a ring and \(S\) a multiplicative subset of \(R\). Then \(R\) is called a \textbf{\(u\)-\(S\)-coherent ring} (with respect to \(s\)) if there exists an element \(s \in S\) such that every finitely generated ideal of \(R\) is \(u\)-\(S\)-finitely presented with respect to \(s\).
	\end{definition}
	
	The relationship between \(u\)-\(S\)-coherent modules and \(u\)-\(S\)-coherent rings is analogous to the classical case.
	
	\begin{corollary}\cite[Corollary 3.4]{Zhang2025b}
		\label{cor:uscoherentsubmod}
		Let \(R\) be a ring and \(S\) a multiplicative subset of \(R\). Let \(N\) be a submodule of an \(R\)-module \(M\). If \(M\) is \(u\)-\(S\)-coherent and \(N\) is \(S\)-finite, then \(M/N\) is \(u\)-\(S\)-coherent.
	\end{corollary}
	
	A useful ideal-theoretic characterization of \(u\)-\(S\)-coherent rings exists, analogous to the classical one for coherent rings.
	
	\begin{proposition}[Zhang \cite{Zhang2025b}, Proposition 3.11]
		\label{prop:uscoherentidealchar}
		Let \(R\) be a ring and \(S\) a multiplicative subset of \(R\). Then the following statements are equivalent:
		\begin{enumerate}
			\item \(R\) is a \(u\)-\(S\)-coherent ring;
			\item there is \(s \in S\) such that \((0:_R r)\) is \(S\)-finite with respect to \(s\) for any \(r \in R\) and the intersection of two finitely generated ideals of \(R\) is \(S\)-finite with respect to \(s\);
			\item there is \(s \in S\) such that \((I:_R b)\) is \(S\)-finite with respect to \(s\) for any element \(b \in R\) and any finitely generated ideal \(I\) of \(R\).
		\end{enumerate}
	\end{proposition}
	
	The uniform property is stronger than its non-uniform counterparts. A \(u\)-\(S\)-coherent ring is both \(S\)-coherent and \(c\)-\(S\)-coherent.
	
	\begin{proposition}[Zhang \cite{Zhang2025b}, Proposition 3.12]
		\label{prop:uscoherentvsothers}
		Let \(R\) be a ring, \(S\) be a multiplicative subset of \(R\). If \(R\) is a \(u\)-\(S\)-coherent ring, then \(R\) is both \(S\)-coherent and \(c\)-\(S\)-coherent.
	\end{proposition}
	
	When \(S\) is a finite set, all these notions coincide.
	
	\begin{proposition}[Zhang \cite{Zhang2025b}, Proposition 3.13]
		\label{prop:uscoherentfiniteS}
		Let \(R\) be a ring and \(S\) a multiplicative subset of \(R\) consisting of finite elements. Then the following statements are equivalent:
		\begin{enumerate}
			\item \(R\) is a \(u\)-\(S\)-coherent ring;
			\item \(R\) is an \(S\)-coherent ring;
			\item \(R\) is a \(c\)-\(S\)-coherent ring.
		\end{enumerate}
	\end{proposition}
	
	Localizing a \(u\)-\(S\)-coherent ring at the uniform element yields a classical coherent ring.
	
	\begin{proposition}\cite[Proposition 3.14]{Zhang2025b}
		\label{prop:uscoherentloc}
		Let \(R\) be a ring and \(S\) a multiplicative subset of \(R\). If \(R\) is a \(u\)-\(S\)-coherent ring with respect to some \(s \in S\), then \(R_s\) is a coherent ring.
	\end{proposition}
	
	The following example, adapted from \cite{Zhang2025b}, demonstrates that a ring can be both \(S\)-coherent and \(c\)-\(S\)-coherent without being \(u\)-\(S\)-coherent, highlighting the strength of the uniform condition.
	
	\begin{example}\cite[Example 3.15]{Zhang2025b}
		\label{ex:coherentnotuscoherent}
		Let \(R\) be a domain such that \(R_s\) is not coherent for any non-zero \(s \in R\). For example, \(R = \mathbb{Q} + x\mathbb{R}[[x]]\). Set \(S = R \setminus \{0\}\). Then \(R\) is both \(S\)-coherent and \(c\)-\(S\)-coherent, but it is not \(u\)-\(S\)-coherent.
	\end{example}

	\subsection{Uniformly S-Artinian Rings and Modules}

Artinian rings are those satisfying the descending chain condition on ideals. The uniform \(S\)-analogue requires a uniform element to control the chain condition.

\begin{definition} \cite[Definition 2.1]{Zhang2025a}
	\label{def:usartinianmodule}
	Let \(R\) be a ring and \(S\) a multiplicative subset of \(R\). An \(R\)-module \(M\) is called a \textbf{\(u\)-\(S\)-Artinian} module (with respect to \(s\)) provided that there exists \(s \in S\) such that each descending chain \(\{M_i\}_{i \in \mathbb{Z}^+}\) of submodules of \(M\) is \(S\)-stationary with respect to \(s\), i.e., there exists \(k \ge 1\) such that \(sM_k \subseteq M_n\) for all \(n \ge k\).
\end{definition}

A ring is \(u\)-\(S\)-Artinian if it is \(u\)-\(S\)-Artinian as a module over itself.

\begin{definition}[Zhang and Qi \cite{Zhang2025a}, Definition 3.1]
	\label{def:usartinianring}
	Let \(R\) be a ring and \(S\) a multiplicative subset of \(R\). Then \(R\) is called a \textbf{\(u\)-\(S\)-Artinian} ring (with respect to \(s\)) provided that \(R\) is a \(u\)-\(S\)-Artinian \(R\)-module.
\end{definition}

The property is invariant under taking the saturation of \(S\).

\begin{proposition} \cite[Proposition 2.3]{Zhang2025a}
	\label{prop:usartiniansat}
	Let \(R\) be a ring, \(S\) a multiplicative subset of \(R\) and \(M\) an \(R\)-module. Let \(S^{*}\) be the saturation of \(S\). Then \(M\) is a \(u\)-\(S\)-Artinian \(R\)-module if and only if \(M\) is a \(u\)-\(S^{*}\)-Artinian \(R\)-module.
\end{proposition}

Localizing a \(u\)-\(S\)-Artinian module at the uniform element yields a classically Artinian module.

\begin{lemma}\cite[Lemma 2.4]{Zhang2025a}
	\label{lem:usartinianloc}
	Let \(R\) be a ring, \(S\) a multiplicative subset of \(R\) and \(M\) an \(R\)-module. If \(M\) is a \(u\)-\(S\)-Artinian \(R\)-module, then there exists an element \(s \in S\) such that \(M_s\) is an Artinian \(R_s\)-module.
\end{lemma}

When \(S\) satisfies the maximal multiple condition, \(S\)-Artinian and \(u\)-\(S\)-Artinian coincide.

\begin{proposition} \cite[Proposition 2.6]{Zhang2025a}
	\label{prop:usartinianiffartinianmmc}
	Let \(R\) be a ring, \(S\) a multiplicative subset of \(R\) satisfying the maximal multiple condition, and \(M\) an \(R\)-module. Then \(M\) is a \(u\)-\(S\)-Artinian module if and only if \(M\) is an \(S\)-Artinian module.
\end{proposition}

However, in general, an \(S\)-Artinian module need not be \(u\)-\(S\)-Artinian, as shown by a valuation domain example.

\begin{example} \cite[Example 2.7]{Zhang2025a}
	\label{ex:sartiniannotusartinian}
	Let \(R\) be a valuation domain whose valuation group is the Hahn product \(G = \prod_{\mathbb{R}} \mathbb{Z}\). Let \(S = R \setminus \{0\}\). Then \(R\) itself is an \(S\)-Artinian \(R\)-module but not \(u\)-\(S\)-Artinian.
\end{example}

The property is invariant under \(u\)-\(S\)-isomorphisms and behaves well in \(u\)-\(S\)-exact sequences.

\begin{lemma} \cite[Lemma 2.8]{Zhang2025a}
	\label{lem:usartinianisoinv}
	Let \(R\) be a ring and \(S\) a multiplicative subset of \(R\). Let \(M\) and \(N\) be \(R\)-modules. If \(M\) is \(u\)-\(S\)-isomorphic to \(N\), then \(M\) is \(u\)-\(S\)-Artinian if and only if \(N\) is \(u\)-\(S\)-Artinian.
\end{lemma}

\begin{proposition} \cite[Proposition 2.9]{Zhang2025a}
	\label{prop:usartinianexact}
	Let \(R\) be a ring and \(S\) a multiplicative subset of \(R\). Let \(0 \to A \to B \to C \to 0\) be an \(S\)-exact sequence. Then \(B\) is \(u\)-\(S\)-Artinian if and only if \(A\) and \(C\) are \(u\)-\(S\)-Artinian. Consequently, a finite direct sum \(\bigoplus_{i=1}^n M_i\) is \(u\)-\(S\)-Artinian if and only if each \(M_i\) is \(u\)-\(S\)-Artinian.
\end{proposition}

A local characterization of classical Artinian modules is obtained.

\begin{proposition} \cite[Proposition 2.10]{Zhang2025a}
	\label{prop:artinianlocalus}
	Let \(R\) be a ring and \(M\) an \(R\)-module. Then the following statements are equivalent:
	\begin{enumerate}
		\item \(M\) is Artinian;
		\item \(M\) is \(u\)-\(\mathfrak{p}\)-Artinian for any \(\mathfrak{p} \in \Spec(R)\);
		\item \(M\) is \(u\)-\(\mathfrak{m}\)-Artinian for any \(\mathfrak{m} \in \Max(R)\).
	\end{enumerate}
\end{proposition}

The notion of a \(u\)-\(S\)-cofinite module is the uniform version of finitely cogenerated modules.

\begin{definition} \cite[Definition 2.11]{Zhang2025a}
	\label{def:uscofinite}
	Let \(R\) be a ring and \(S\) a multiplicative subset of \(R\). An \(R\)-module \(M\) is called \textbf{\(u\)-\(S\)-cofinite} (with respect to \(s\)) if there is an \(s \in S\) such that for each nonempty family of submodules \(\{M_i\}_{i \in \Delta}\) of \(M\), \(\bigcap_{i \in \Delta} M_i = 0\) implies that \(s(\bigcap_{i \in \Delta'} M_i) = 0\) for a finite subset \(\Delta' \subseteq \Delta\).
\end{definition}

The following theorem gives several characterizations of \(u\)-\(S\)-Artinian modules, including conditions involving minimal elements and cofiniteness of quotients.

\begin{theorem}\cite[Theorem 2.15]{Zhang2025a}
	\label{thm:usartinianchar}
	Let \(R\) be a ring, \(S\) a multiplicative subset of \(R\) and \(M\) an \(R\)-module. Let \(s \in S\). Then the following statements are equivalent:
	\begin{enumerate}
		\item \(M\) is a \(u\)-\(S\)-Artinian module with respect to \(s\);
		\item \(M\) satisfies the (S-MIN)-condition with respect to \(s\) (every nonempty family of submodules has an \(S\)-minimal element w.r.t. \(s\));
		\item For any nonempty family \(\{N_i\}_{i \in \Gamma}\) of submodules of \(M\), there is a finite subset \(\Gamma_0 \subseteq \Gamma\) such that \(s \bigcap_{i \in \Gamma_0} N_i \subseteq \bigcap_{i \in \Gamma} N_i\);
		\item Every factor module \(M/N\) is \(u\)-\(S\)-cofinite with respect to \(s\).
	\end{enumerate}
\end{theorem}

For a ring with a regular multiplicative set, \(u\)-\(S\)-Artinian implies classically Artinian.

\begin{proposition} \cite[Proposition 3.2]{Zhang2025a}
	\label{prop:usartinianregularimpliesart}
	Let \(R\) be a ring and \(S\) a regular multiplicative subset of \(R\). If \(R\) is a \(u\)-\(S\)-Artinian ring, then \(R\) is an Artinian ring.
\end{proposition}

The property is well-behaved under finite direct products.

\begin{proposition} \cite[Proposition 3.3]{Zhang2025a}
	\label{prop:usartinianproduct}
	Let \(R = R_1 \times R_2\) be direct product of rings \(R_1\) and \(R_2\), and \(S = S_1 \times S_2\) a direct product of multiplicative subsets. Then \(R\) is a \(u\)-\(S\)-Artinian ring if and only if \(R_i\) is a \(u\)-\(S_i\)-Artinian ring for each \(i = 1,2\).
\end{proposition}

An example of a non-Artinian \(u\)-\(S\)-Artinian ring is constructed as a direct product.

\begin{example}\cite[Example 3.4]{Zhang2025a}
	\label{ex:usartiniannotart}
	Let \(R = R_1 \times R_2\), where \(R_1\) is Artinian and \(R_2\) is not Artinian. Set \(S = \{1\} \times \{1,0\}\). Then \(R\) is \(u\)-\(S\)-Artinian but not Artinian.
\end{example}

A key result is the structure theorem for \(u\)-\(S\)-Artinian rings, which parallels the classical theorem: a ring is \(u\)-\(S\)-Artinian if and only if it is \(u\)-\(S\)-Noetherian, its \(u\)-\(S\)-Jacobson radical \(\Jac_S(R)\) is \(S\)-nilpotent, and the quotient by this radical is \(u\)-\(S\)-semisimple.

First, we need the definition of the \(u\)-\(S\)-Jacobson radical.

\begin{definition} \cite[Definition 4.1 and 4.2]{Zhang2025a}
	\label{def:usjacobson}
	Let \(R\) be a ring, \(S\) a multiplicative subset of \(R\). For an \(R\)-module \(M\), the \textbf{\(u\)-\(S\)-Jacobson radical} \(\Jac_S(M)\) is defined as the intersection of all \(u\)-\(S\)-maximal submodules of \(M\). For the ring \(R\), \(\Jac_S(R)\) is the \(u\)-\(S\)-Jacobson radical of \(R\) as an \(R\)-module.
\end{definition}

A key lemma shows that \(u\)-\(S\)-Artinian modules have a finite decomposition into \(u\)-\(S\)-simple modules, up to a \(u\)-\(S\)-torsion submodule.

\begin{proposition} \cite[Proposition 4.4]{Zhang2025a}
	\label{prop:usartiniandecomp}
	Suppose \(M\) is a \(u\)-\(S\)-Artinian \(R\)-module with \(\Jac_S(M)\) \(u\)-\(S\)-torsion. Then there exists \(T \subseteq M\) such that \(sT = 0\) and \(M/T \subseteq \bigoplus_{i=1}^{t} S_i\) where each \(S_i\) is \(u\)-\(S\)-simple with respect to \(s\) for some \(s \in S\). Consequently, \(M\) is a \(u\)-\(S\)-Noetherian \(R\)-module.
\end{proposition}

From this, we deduce that for a \(u\)-\(S\)-Artinian ring, the quotient by its radical is \(u\)-\(S\)-semisimple.

\begin{proposition} \cite[Proposition 4.5]{Zhang2025a}
	\label{prop:usartinianquotss}
	Suppose \(R\) is a \(u\)-\(S\)-Artinian ring. Then \(R/\Jac_S(R)\) is a \(u\)-\(S/\Jac_S(R)\)-semisimple ring.
\end{proposition}

The \(u\)-\(S\)-Jacobson radical of a \(u\)-\(S\)-Artinian ring is \(S\)-nilpotent.

\begin{proposition} \cite[Proposition 4.8]{Zhang2025a}
	\label{prop:usartinianradnilp}
	Let \(R\) be a ring and \(S\) a multiplicative subset of \(R\). Suppose \(R\) is a \(u\)-\(S\)-Artinian ring. Then \(\Jac_S(R)\) is \(S\)-nilpotent.
\end{proposition}

Finally, the main structure theorem.

\begin{theorem} \cite[Theorem 4.9]{Zhang2025a}
	\label{thm:usartinianstruct}
	Let \(R\) be a ring and \(S\) a multiplicative subset of \(R\). Then the following statements are equivalent:
	\begin{enumerate}
		\item \(R\) is a \(u\)-\(S\)-Artinian ring;
		\item \(R\) is a \(u\)-\(S\)-Noetherian ring, \(\Jac_S(R)\) is \(S\)-nilpotent, and \(R/\Jac_S(R)\) is a \(u\)-\(S/\Jac_S(R)\)-semisimple ring.
	\end{enumerate}
\end{theorem}

The property is also studied for idealizations, showing that the \(u\)-\(S\)-Artinian property of \(R(+)M\) is equivalent to \(R\) being \(u\)-\(S\)-Artinian and \(M\) being \(S\)-finite.

\begin{proposition} \cite[Proposition 4.10]{Zhang2025a}
	\label{prop:usartinianidealization}
	Let \(R\) be a commutative ring, \(S\) a multiplicative subset of \(R\) and \(M\) an \(R\)-module. Then \(R(+)M\) is a \(u\)-\(S(+)M\)-Artinian ring if and only if \(R\) is a \(u\)-\(S\)-Artinian ring and \(M\) is an \(S\)-finite \(R\)-module.
\end{proposition}
	
	\section{Uniformly S-versions of Classical Modules}
	
	With the foundation of \(u\)-\(S\)-torsion and \(u\)-\(S\)-exact sequences, we can now define and study the uniform analogues of the fundamental module classes in homological algebra. These definitions are direct translations of the classical ones, replacing exactness with \(u\)-\(S\)-exactness.
	
	\subsection{Uniformly S-Flat Modules}
	
	The notion of flatness is central to homological algebra. Its uniform \(S\)-analogue was the first to be systematically developed.
	
	\begin{definition}\cite[Definition 3.1]{Zhang2022}
		\label{def:usflat}
		Let \(R\) be a ring, \(S\) a multiplicative subset of \(R\). An \(R\)-module \(F\) is called \textbf{\(u\)-\(S\)-flat} provided that for any \(u\)-\(S\)-exact sequence \(0 \to A \to B \to C \to 0\), the induced sequence \(0 \to A \otimes_R F \to B \otimes_R F \to C \otimes_R F \to 0\) is \(u\)-\(S\)-exact.
	\end{definition}
	
	The following theorem provides several equivalent characterizations, the most important of which is that \(\Tor_1^R(M,F)\) is \(u\)-\(S\)-torsion for all modules \(M\). This is the perfect \(u\)-\(S\)-analogue of the classical flatness criterion.
	
	\begin{theorem} \cite[Theorem 3.2]{Zhang2022}
		\label{thm:usflatchar}
		Let \(R\) be a ring, \(S\) a multiplicative subset of \(R\) and \(F\) an \(R\)-module. The following statements are equivalent:
		\begin{enumerate}
			\item \(F\) is \(u\)-\(S\)-flat;
			\item for any short exact sequence \(0 \to A \xrightarrow{f} B \xrightarrow{g} C \to 0\), the induced sequence \(0 \to A \otimes_R F \xrightarrow{f \otimes_R F} B \otimes_R F \xrightarrow{g \otimes_R F} C \otimes_R F \to 0\) is \(u\)-\(S\)-exact;
			\item \(\Tor_1^R(M,F)\) is \(u\)-\(S\)-torsion for any \(R\)-module \(M\);
			\item \(\Tor_n^R(M,F)\) is \(u\)-\(S\)-torsion for any \(R\)-module \(M\) and \(n \ge 1\).
		\end{enumerate}
	\end{theorem}
	
	An important example from \cite{Zhang2022} shows that, unlike in the classical case, checking \(\Tor_1\) against cyclic modules is insufficient to guarantee \(u\)-\(S\)-flatness.
	
	\begin{example}\cite[Example 3.3]{Zhang2022}
		\label{ex:torcycleinsufficient}
		Let \(R = \mathbb{Z}\), \(p\) a prime, and \(S = \{p^n \mid n \ge 0\}\). Let \(M = \mathbb{Z}_{(p)}/\mathbb{Z}\). Then \(\Tor_1^R(R/I, M)\) is \(u\)-\(S\)-torsion for any ideal \(I\) of \(R\). However, \(M\) is not \(u\)-\(S\)-flat because \(\Tor_1^R(\mathbb{Q}/\mathbb{Z}, M) \cong M\) is not \(u\)-\(S\)-torsion.
	\end{example}
	
	The class of \(u\)-\(S\)-flat modules has several nice closure properties.
	
	\begin{proposition} \cite[Proposition 3.4]{Zhang2022}
		\label{prop:usflatclosure}
		Let \(R\) be a ring and \(S\) a multiplicative subset of \(R\). Then the following statements hold.
		\begin{enumerate}
			\item Any pure quotient of \(u\)-\(S\)-flat modules is \(u\)-\(S\)-flat.
			\item Any finite direct sum of \(u\)-\(S\)-flat modules is \(u\)-\(S\)-flat.
			\item Let \(0 \to A \xrightarrow{f} B \xrightarrow{g} C \to 0\) be a \(u\)-\(S\)-exact sequence. If \(A\) and \(C\) are \(u\)-\(S\)-flat modules, so is \(B\).
			\item Let \(A \to B\) be a \(u\)-\(S\)-isomorphism. If one of \(A\) and \(B\) is \(u\)-\(S\)-flat, so is the other.
			\item Let \(0 \to A \xrightarrow{f} B \xrightarrow{g} C \to 0\) be a \(u\)-\(S\)-exact sequence. If \(B\) and \(C\) are \(u\)-\(S\)-flat, then \(A\) is \(u\)-\(S\)-flat.
		\end{enumerate}
	\end{proposition}
	
	Unlike classical flatness, the \(u\)-\(S\)-flat property is not preserved under arbitrary direct sums or direct limits, as shown by a counterexample from \cite{Zhang2022}.
	
	\begin{remark}\cite[Remark 3.5]{Zhang2022}
		\label{rem:usflatnotclosed}
		Let \(R = \mathbb{Z}\), \(p\) a prime, and \(S = \{p^n \mid n \ge 0\}\). The modules \(F_n = \mathbb{Z}/\langle p^n \rangle\) are \(u\)-\(S\)-torsion, hence \(u\)-\(S\)-flat. However, their direct limit \(\varinjlim F_n \cong \mathbb{Z}_{(p)}/\mathbb{Z}\) is not \(u\)-\(S\)-flat. Similarly, the direct sum \(\bigoplus_{n=1}^{\infty} \mathbb{Z}/\langle p^n \rangle\) is not \(u\)-\(S\)-flat.
	\end{remark}
	
	Localizing a \(u\)-\(S\)-flat module at \(S\) yields a classically flat module.
	
	\begin{corollary}\cite[Corollary 3.6]{Zhang2022}
		\label{cor:usflatloc}
		Let \(R\) be a ring and \(S\) a multiplicative subset of \(R\). If \(F\) is \(u\)-\(S\)-flat over \(R\), then \(F_S\) is flat over \(R_S\).
	\end{corollary}
	
	The converse is false in general, as shown by the next remark.
	
	\begin{remark} \cite[Remark 3.7]{Zhang2022}
		\label{rem:flatlocnotusflat}
		Consider \(R=\mathbb{Z}\), \(S=\{p^n \mid n\ge 0\}\), and \(M = \mathbb{Z}_{(p)}/\mathbb{Z}\). Then \(M_S = 0\) is flat over \(\mathbb{Z}_S\), but \(M\) is not \(u\)-\(S\)-flat over \(\mathbb{Z}\).
	\end{remark}
	
	However, when \(S\) is a finite set, the two properties are equivalent.
	
	\begin{proposition} \cite[Proposition 3.8]{Zhang2022}
		\label{prop:usflatlocfinite}
		Let \(R\) be a ring and \(F\) an \(R\)-module. Let \(S\) be a multiplicative subset of \(R\) consisting of finite elements. Then \(F\) is \(u\)-\(S\)-flat over \(R\) if and only if \(F_S\) is flat over \(R_S\).
	\end{proposition}
	
	A local characterization of classical flatness emerges from the uniform theory.
	
	\begin{proposition}\cite[Proposition 3.9]{Zhang2022}
		\label{prop:flatlocalus}
		Let \(R\) be a ring and \(F\) an \(R\)-module. Then the following statements are equivalent:
		\begin{enumerate}
			\item \(F\) is flat;
			\item \(F\) is \(u\)-\(\mathfrak{p}\)-flat for any \(\mathfrak{p} \in \Spec(R)\);
			\item \(F\) is \(u\)-\(\mathfrak{m}\)-flat for any \(\mathfrak{m} \in \Max(R)\).
		\end{enumerate}
	\end{proposition}
	
	Further characterizations of \(u\)-\(S\)-flatness, particularly useful when \(S\) is regular, are given in \cite{Zhang2025b}.
	
	\begin{proposition} \cite[Proposition 4.3]{Zhang2025b}
		\label{prop:usflatmorechar}
		Let \(R\) be a ring and \(S\) be a multiplicative subset of \(R\). Then the following statements are equivalent:
		\begin{enumerate}
			\item \(F\) is \(u\)-\(S\)-flat;
			\item there exists an element \(s\in S\) such that \(\Tor_1^R(N,F)\) is uniformly \(S\)-torsion with respect to \(s\) for any finitely presented \(R\)-module \(N\);
			\item \(\Hom_R(F,E)\) is \(u\)-\(S\)-injective for any injective module \(E\);
			\item \(\Hom_R(F,E)\) is \(u\)-\(S\)-absolutely pure for any injective module \(E\);
			\item if \(E\) is an injective cogenerator, then \(\Hom_R(F,E)\) is \(u\)-\(S\)-injective;
			\item if \(E\) is an injective cogenerator, then \(\Hom_R(F,E)\) is \(u\)-\(S\)-absolutely pure.
		\end{enumerate}
		Moreover, if \(S\) is regular, then all above are equivalent to:
		\begin{enumerate}
			\item[7.] there exists \(s \in S\) such that \(\Tor_1^R(R/I, F)\) is uniformly \(S\)-torsion with respect to \(s\) for any ideal \(I\) of \(R\);
			\item[8.] there exists \(s \in S\) such that, for any ideal \(I\) of \(R\), the natural homomorphism \(\sigma_I: I \otimes_R F \to IF\) is a \(u\)-\(S\)-isomorphism with respect to \(s\);
			\item[9.] there exists \(s \in S\) such that \(\Tor_1^R(R/K, F)\) is uniformly \(S\)-torsion with respect to \(s\) for any finitely generated ideal \(K\) of \(R\);
			\item[10.] there exists \(s \in S\) such that, for any finitely generated ideal \(K\) of \(R\), the natural homomorphism \(\sigma_K: K \otimes_R F \to KF\) is a \(u\)-\(S\)-isomorphism with respect to \(s\).
		\end{enumerate}
	\end{proposition}
	
	\subsection{Uniformly S-Projective Modules}
	
	Dual to the notion of flatness is projectivity. The uniform \(S\)-analogue is defined by the behavior of the \(\Hom\) functor.
	
	\begin{definition} \cite[Definition 2.8]{Zhang2023a}
		\label{def:usprojective}
		Let \(R\) be a ring and \(S\) a multiplicative subset of \(R\). An \(R\)-module \(P\) is called \textbf{\(u\)-\(S\)-projective} provided that the induced sequence
		\[
		0 \to \Hom_R(P,A) \to \Hom_R(P,B) \to \Hom_R(P,C) \to 0
		\]
		is \(u\)-\(S\)-exact for any \(u\)-\(S\)-exact sequence \(0 \to A \to B \to C \to 0\).
	\end{definition}
	
	The following theorem provides the key characterizations, parallel to those for \(u\)-\(S\)-injective modules. The equivalence with \(\Ext_R^1(P,M)\) being \(u\)-\(S\)-torsion is particularly useful.
	
	\begin{theorem} \cite[Theorem 2.9]{Zhang2023a}
		\label{thm:usprojectivechar}
		Let \(R\) be a ring, \(S\) a multiplicative subset of \(R\) and \(P\) an \(R\)-module. Then the following statements are equivalent:
		\begin{enumerate}
			\item \(P\) is \(u\)-\(S\)-projective;
			\item for any \(u\)-\(S\)-epimorphism \(B \xrightarrow{g} C\) there exists \(s \in S\) such that for any \(R\)-homomorphism \(h:P \to C\), there exists an \(R\)-homomorphism \(\alpha :P \to B\) satisfying \(sh = g \circ \alpha\);
			\item for any short exact sequence \(0 \to A \xrightarrow{f} B \xrightarrow{g} C \to 0\), the induced sequence \(0 \to \Hom_R(P,A) \xrightarrow{f_*} \Hom_R(P,B) \xrightarrow{g_*} \Hom_R(P,C) \to 0\) is \(u\)-\(S\)-exact;
			\item \(\Ext_R^1(P,M)\) is \(u\)-\(S\)-torsion for any \(R\)-module \(M\);
			\item \(\Ext_R^n(P,M)\) is \(u\)-\(S\)-torsion for any \(R\)-module \(M\) and \(n \ge 1\).
		\end{enumerate}
	\end{theorem}
	
	The relation between \(u\)-\(S\)-projective and \(u\)-\(S\)-split exact sequences is clarified in the following corollary.
	
	\begin{corollary} \cite[Corollary 2.10]{Zhang2023a}
		\label{cor:usprojectivesplit}
		Let \(R\) be a ring, \(S\) a multiplicative subset of \(R\) and \(P\) an \(R\)-module. Then the following statements hold:
		\begin{enumerate}
			\item If \(P\) is \(u\)-\(S\)-projective, then any \(u\)-\(S\)-short exact sequence \(0 \to A \xrightarrow{f} B \xrightarrow{g} P \to 0\) is \(u\)-\(S\)-split;
			\item If there is \(s \in S\) such that any short exact sequence \(0 \to A \xrightarrow{f} B \xrightarrow{g} P \to 0\) is \(u\)-\(S\)-split with respect to \(s\), then \(P\) is \(u\)-\(S\)-projective.
		\end{enumerate}
	\end{corollary}
	
	Just as with flatness, \(u\)-\(S\)-torsion modules are trivially \(u\)-\(S\)-projective.
	
	\begin{corollary} \cite[Corollary 2.11]{Zhang2023a}
		\label{cor:ustorsionproj}
		Let \(R\) be a ring and \(S\) a multiplicative subset of \(R\). Let \(P\) be a \(u\)-\(S\)-torsion \(R\)-module or a projective \(R\)-module. Then \(P\) is \(u\)-\(S\)-projective.
	\end{corollary}
	
	The \(u\)-\(S\)-projective property is also preserved under finite direct sums and behaves well in \(u\)-\(S\)-exact sequences.
	
	\begin{proposition} \cite[Proposition 2.14]{Zhang2023a}
		\label{prop:usprojectiveclosure}
		Let \(R\) be a ring and \(S\) a multiplicative subset of \(R\). Then the following statements hold.
		\begin{enumerate}
			\item Any finite direct sum of \(u\)-\(S\)-projective modules is \(u\)-\(S\)-projective.
			\item Let \(0 \to A \xrightarrow{f} B \xrightarrow{g} C \to 0\) be a \(u\)-\(S\)-exact sequence. If \(C\) is \(u\)-\(S\)-projective, then \(A\) is \(u\)-\(S\)-projective if and only if so is \(B\).
			\item Let \(A \to B\) be a \(u\)-\(S\)-isomorphism. Then \(A\) is \(u\)-\(S\)-projective if and only if \(B\) is \(u\)-\(S\)-projective.
			\item Let \(0 \to A \xrightarrow{f} B \xrightarrow{g} C \to 0\) be a \(u\)-\(S\)-split \(u\)-\(S\)-exact sequence. If \(B\) is \(u\)-\(S\)-projective, then \(A\) and \(C\) are \(u\)-\(S\)-projective.
		\end{enumerate}
	\end{proposition}
	
	However, arbitrary direct sums of \(u\)-\(S\)-projective modules need not be \(u\)-\(S\)-projective.
	
	\begin{example}\cite[Example 2.15]{Zhang2023a}
		\label{ex:usprojectivesumn}
		Let \(R = \mathbb{Z}\), \(S = \{p^n \mid n \in \mathbb{N}\}\). The modules \(M_n = \mathbb{Z}/\langle p^n \rangle\) are \(u\)-\(S\)-torsion, hence \(u\)-\(S\)-projective. Their direct sum \(N = \bigoplus_{n=1}^{\infty} M_n\) is not \(u\)-\(S\)-projective because \(\Ext_R^1(N,N)\) is not \(u\)-\(S\)-torsion.
	\end{example}
	
	The relationship between \(u\)-\(S\)-projective and \(u\)-\(S\)-flat modules is as expected.
	
	\begin{proposition} \cite[Proposition 2.13]{Zhang2023a}
		\label{prop:usprojectiveimpliesflat}
		Let \(R\) be a ring and \(S\) a multiplicative subset of \(R\). If \(P\) is a \(u\)-\(S\)-projective \(R\)-module, then \(P\) is \(u\)-\(S\)-flat.
	\end{proposition}
	
	A local characterization of classical projective modules is given next.
	
	\begin{proposition} \cite[Proposition 2.16]{Zhang2023a}
		\label{prop:projectivelocalus}
		Let \(R\) be a ring and \(P\) an \(R\)-module. Then the following statements are equivalent:
		\begin{enumerate}
			\item \(P\) is projective;
			\item \(P\) is \(u\)-\(\mathfrak{p}\)-projective for any \(\mathfrak{p} \in \Spec(R)\);
			\item \(P\) is \(u\)-\(\mathfrak{m}\)-projective for any \(\mathfrak{m} \in \Max(R)\).
		\end{enumerate}
	\end{proposition}
	A connection between \(u\)-\(S\)-projective, \(u\)-\(S\)-flat, and \(u\)-\(S\)-finitely presented modules is given in the next proposition.

\begin{proposition}\cite[Proposition 2.8]{Zhang2025b}
	\label{prop:projflatfp}
	Let \(R\) be a ring and \(S\) be a multiplicative subset of \(R\). Then the following statements hold.
	\begin{enumerate}
		\item Every \(S\)-finite \(u\)-\(S\)-projective module is \(u\)-\(S\)-finitely presented.
		\item Every \(u\)-\(S\)-finitely presented \(u\)-\(S\)-flat module is \(u\)-\(S\)-projective.
	\end{enumerate}
\end{proposition}
	
	\subsection{Uniformly S-Injective Modules}
	
	The dual notion of projectivity is injectivity. The uniform \(S\)-analogue of injective modules was first studied in the context of \(u\)-\(S\)-Noetherian rings.
	
	\begin{definition}\cite[Definition 4.1]{Qi2022}
		\label{def:usinjective}
		Let \(R\) be a ring and \(S\) a multiplicative subset of \(R\). An \(R\)-module \(E\) is said to be \textbf{\(u\)-\(S\)-injective} provided that the induced sequence
		\[
		0 \to \Hom_R(C,E) \to \Hom_R(B,E) \to \Hom_R(A,E) \to 0
		\]
		is \(u\)-\(S\)-exact for any \(u\)-\(S\)-exact sequence \(0 \to A \to B \to C \to 0\) of \(R\)-modules.
	\end{definition}
	
	A key technical lemma is that \(\Ext\) with a \(u\)-\(S\)-torsion module is again \(u\)-\(S\)-torsion.
	
	\begin{lemma}\cite[Lemma 4.2]{Qi2022}
		\label{lem:extustorsion}
		Let \(R\) be a ring and \(S\) a multiplicative subset of \(R\). If \(T\) is a uniformly \(S\)-torsion module, then \(\Ext_R^n(T,M)\) and \(\Ext_R^n(M,T)\) are uniformly \(S\)-torsion for any \(R\)-module \(M\) and any \(n \ge 0\).
	\end{lemma}
	
	The characterization of \(u\)-\(S\)-injective modules mirrors the classical one.
	
	\begin{theorem}\cite[Theorem 4.3]{Qi2022}
		\label{thm:usinjectivechar}
		Let \(R\) be a ring, \(S\) a multiplicative subset of \(R\), and \(E\) an \(R\)-module. Then the following conditions are equivalent:
		\begin{enumerate}
			\item \(E\) is \(u\)-\(S\)-injective;
			\item for any short exact sequence \(0 \to A \xrightarrow{f} B \xrightarrow{g} C \to 0\) of \(R\)-modules, the induced sequence \(0 \to \Hom_R(C,E) \xrightarrow{g^*} \Hom_R(B,E) \xrightarrow{f^*} \Hom_R(A,E) \to 0\) is \(u\)-\(S\)-exact;
			\item \(\Ext_R^1(M,E)\) is uniformly \(S\)-torsion for any \(R\)-module \(M\);
			\item \(\Ext_R^n(M,E)\) is uniformly \(S\)-torsion for any \(R\)-module \(M\) and \(n \ge 1\).
		\end{enumerate}
	\end{theorem}
	
	A consequence of this theorem is that both \(u\)-\(S\)-torsion and injective modules are \(u\)-\(S\)-injective.

	\begin{corollary} \cite[Corollary 4.4]{Qi2022}
		\label{cor:ustorsionandinjareusinj}
		Let \(R\) be a ring and \(S\) a multiplicative subset of \(R\). Suppose that \(E\) is a uniformly \(S\)-torsion \(R\)-module or an injective \(R\)-module. Then \(E\) is \(u\)-\(S\)-injective.
	\end{corollary}
	
	Similar to the flat and projective cases, checking injectivity against cyclic modules is insufficient, as shown by an example involving \(p\)-adic integers.
	
	\begin{example}\cite[Example 4.5]{Qi2022}
		\label{ex:extcycleinsufficient}
		Let \(R = \mathbb{Z}\), \(p\) a prime, and \(S = \{p^n \mid n \ge 0\}\). Let \(J_p\) be the additive group of \(p\)-adic integers. Then \(\Ext_R^1(R/I, J_p)\) is \(u\)-\(S\)-torsion for any ideal \(I\) of \(R\). However, \(J_p\) is not \(u\)-\(S\)-injective because \(\Ext_R^1(\mathbb{Z}(p^{\infty}), J_p) \cong J_p\) is not \(u\)-\(S\)-torsion.
	\end{example}
	
	The class of \(u\)-\(S\)-injective modules has closure properties analogous to those of \(u\)-\(S\)-projective modules.
	
	\begin{proposition} \cite[Proposition 4.7]{Qi2022}
		\label{prop:usinjectiveclosure}
		Let \(R\) be a ring and \(S\) a multiplicative subset of \(R\). Then the following statements hold.
		\begin{enumerate}
			\item Any finite direct sum of \(u\)-\(S\)-injective modules is \(u\)-\(S\)-injective.
			\item Let \(0 \to A \xrightarrow{f} B \xrightarrow{g} C \to 0\) be a \(u\)-\(S\)-exact sequence of \(R\)-modules. If \(A\) and \(C\) are \(u\)-\(S\)-injective modules, so is \(B\).
			\item Let \(A \to B\) be a \(u\)-\(S\)-isomorphism of \(R\)-modules. If one of \(A\) and \(B\) is \(u\)-\(S\)-injective, so is the other.
			\item Let \(0 \to A \xrightarrow{f} B \xrightarrow{g} C \to 0\) be a \(u\)-\(S\)-exact sequence of \(R\)-modules. If \(A\) and \(B\) are \(u\)-\(S\)-injective, then \(C\) is \(u\)-\(S\)-injective.
		\end{enumerate}
	\end{proposition}
	
	A local characterization of classical injective modules is also available.
	
	\begin{proposition} \cite[Proposition 4.8]{Qi2022}
		\label{prop:injectivelocalus}
		Let \(R\) be a ring and \(E\) an \(R\)-module. Then the following conditions are equivalent:
		\begin{enumerate}
			\item \(E\) is injective;
			\item \(E\) is \(\mathfrak{p}\)-injective for any \(\mathfrak{p} \in \Spec(R)\);
			\item \(E\) is \(\mathfrak{m}\)-injective for any \(\mathfrak{m} \in \Max(R)\).
		\end{enumerate}
	\end{proposition}
	
	A uniform version of Baer's Criterion for injectivity is given next. It requires the module to be \(S\)-divisible.
	
	\begin{proposition} \cite[Proposition 4.9]{Qi2022}
		\label{prop:usinjectivebaer}
		(Baer's Criterion for \(u\)-\(S\)-injective modules)
		Let \(R\) be a ring, \(S\) a multiplicative subset of \(R\), and \(E\) an \(R\)-module. If \(E\) is a \(u\)-\(S\)-injective module, then there exists \(s \in S\) such that \(s \Ext_R^1(R/I, E) = 0\) for any ideal \(I\) of \(R\). Moreover, if \(E\) is \(S\)-divisible, then the converse also holds.
	\end{proposition}
	
	The culmination of these properties is the Cartan-Eilenberg-Bass theorem for \(u\)-\(S\)-Noetherian rings, which characterizes them by the behavior of direct sums and direct unions of injective modules.
	
	\begin{theorem}\cite{Qi2022}
		\label{thm:usnoetherianbass}
		(Cartan-Eilenberg-Bass Theorem for uniformly \(S\)-Noetherian rings)
		Let \(R\) be a ring and \(S\) a regular multiplicative subset of \(R\). Then the following conditions are equivalent:
		\begin{enumerate}
			\item \(R\) is uniformly \(S\)-Noetherian;
			\item any direct sum of injective modules is \(u\)-\(S\)-injective;
			\item any direct union of injective modules is \(u\)-\(S\)-injective.
		\end{enumerate}
	\end{theorem}
	
	\subsection{Uniformly S-Absolutely Pure Modules}
	
	Absolutely pure (or FP-injective) modules are a generalization of injective modules. They are defined by the vanishing of \(\Ext_R^1(N, -)\) for all finitely presented modules \(N\). The uniform \(S\)-analogue relaxes this condition.
	
	\begin{definition}\cite[Definition 3.1]{Zhang2023d}
		\label{def:usabsolutelypure}
		Let \(R\) be a ring and \(S\) a multiplicative subset of \(R\). An \(R\)-module \(E\) is said to be \textbf{\(u\)-\(S\)-absolutely pure} provided that any short \(u\)-\(S\)-exact sequence \(0 \to E \to B \to C \to 0\) beginning with \(E\) is \(u\)-\(S\)-pure.
	\end{definition}
	
	The following theorem provides several characterizations, including the key homological one: \(\Ext_R^1(N,E)\) is \(u\)-\(S\)-torsion for all finitely presented modules \(N\).
	
	\begin{theorem}\cite[Theorem 3.2]{Zhang2023d}
		\label{thm:usabsolutelypurechar}
		Let \(R\) be a ring, \(S\) a multiplicative subset of \(R\) and \(E\) an \(R\)-module. Then the following statements are equivalent:
		\begin{enumerate}
			\item \(E\) is \(u\)-\(S\)-absolutely pure;
			\item any short exact sequence \(0 \to E \to B \to C \to 0\) beginning with \(E\) is \(u\)-\(S\)-pure;
			\item \(E\) is a \(u\)-\(S\)-pure submodule in every \(u\)-\(S\)-injective module containing \(E\);
			\item \(E\) is a \(u\)-\(S\)-pure submodule in every injective module containing \(E\);
			\item \(E\) is a \(u\)-\(S\)-pure submodule in its injective envelope;
			\item there exists an element \(s \in S\) satisfying that for any finitely presented \(R\)-module \(N\), \(\Ext_R^1(N,E)\) is \(u\)-\(S\)-torsion with respect to \(s\);
			\item there exists an element \(s \in S\) satisfying that if \(P\) is finitely generated projective, \(K\) is a finitely generated submodule of \(P\) and \(f:K \to E\) is an \(R\)-homomorphism, then there is an \(R\)-homomorphism \(g:P \to E\) such that \(sf = gi\).
		\end{enumerate}
	\end{theorem}
	
	The class of \(u\)-\(S\)-absolutely pure modules is well-behaved under various constructions.
	
	\begin{proposition}[Zhang \cite{Zhang2023d}, Proposition 3.3]
		\label{prop:usapclosure}
		Let \(R\) be a ring and \(S\) a multiplicative subset of \(R\). Then the following statements hold.
		\begin{enumerate}
			\item Any absolutely pure module and any \(u\)-\(S\)-injective module are \(u\)-\(S\)-absolutely pure.
			\item Any finite direct sum of \(u\)-\(S\)-absolutely pure modules is \(u\)-\(S\)-absolutely pure.
			\item Let \(0 \to A \xrightarrow{f} B \xrightarrow{g} C \to 0\) be a \(u\)-\(S\)-exact sequence. If \(A\) and \(C\) are \(u\)-\(S\)-absolutely pure modules, so is \(B\).
			\item The class of \(u\)-\(S\)-absolutely pure modules is closed under \(u\)-\(S\)-isomorphisms.
			\item Let \(0 \to A \to B \to C \to 0\) be a \(u\)-\(S\)-pure \(u\)-\(S\)-exact sequence. If \(B\) is \(u\)-\(S\)-absolutely pure, so is \(A\).
		\end{enumerate}
	\end{proposition}
	
	A local characterization of classical absolutely pure modules is also possible.
	
	\begin{proposition}[Zhang \cite{Zhang2023d}, Proposition 3.4]
		\label{prop:aplocalus}
		Let \(R\) be a ring and \(E\) an \(R\)-module. Then the following statements are equivalent:
		\begin{enumerate}
			\item \(E\) is absolutely pure;
			\item \(E\) is \(u\)-\(\mathfrak{p}\)-absolutely pure for any \(\mathfrak{p} \in \Spec(R)\);
			\item \(E\) is \(u\)-\(\mathfrak{m}\)-absolutely pure for any \(\mathfrak{m} \in \Max(R)\).
		\end{enumerate}
	\end{proposition}
	
	The following theorem establishes a beautiful connection: a ring is \(u\)-\(S\)-von Neumann regular if and only if every module is \(u\)-\(S\)-absolutely pure. This parallels the classical result that a ring is von Neumann regular iff every module is absolutely pure.
	
	\begin{theorem}[Zhang \cite{Zhang2023d}, Theorem 3.5]
		\label{thm:usvnriffallap}
		A ring \(R\) is uniformly \(S\)-von Neumann regular if and only if any \(R\)-module is \(u\)-\(S\)-absolutely pure.
	\end{theorem}
	
	A concrete example shows that \(u\)-\(S\)-absolutely pure modules need not be absolutely pure.
	
	\begin{example}[Zhang \cite{Zhang2023d}, Example 3.6]
		\label{ex:usapnotap}
		Let \(T = \mathbb{Z}_2 \times \mathbb{Z}_2\), \(s = (1,0)\), and \(R = T[x]/\langle sx, x^2 \rangle\) with multiplicative set \(S = \{1, \bar{s}\}\). Then \(R\) is \(u\)-\(S\)-von Neumann regular but not von Neumann regular. Hence, by Theorem \ref{thm:usvnriffallap}, every \(R\)-module is \(u\)-\(S\)-absolutely pure, but not every module is absolutely pure.
	\end{example}
	
	Finally, a characterization of \(u\)-\(S\)-Noetherian rings in terms of \(u\)-\(S\)-absolutely pure modules is given. This is a direct analogue of a classical theorem about Noetherian rings.
	
	\begin{theorem}[Zhang \cite{Zhang2023d}, Theorem 3.7]
		\label{thm:usnoetherianiffapinjective}
		Let \(R\) be a ring, \(S\) a regular multiplicative subset of \(R\). Then the following statements are equivalent:
		\begin{enumerate}
			\item \(R\) is a uniformly \(S\)-Noetherian ring;
			\item any \(u\)-\(S\)-absolutely pure module is \(u\)-\(S\)-injective;
			\item any absolutely pure module is \(u\)-\(S\)-injective.
		\end{enumerate}
	\end{theorem}
	
	Like flat and projective modules, the \(u\)-\(S\)-absolutely pure property is not closed under arbitrary direct sums or products, as shown by a counterexample.
	
	\begin{example}\cite[Example 3.8]{Zhang2023d}
		\label{ex:usapnotclosed}
		Let \(R = \mathbb{Z}\), \(S = \{p^n \mid n \ge 0\}\). The cyclic groups \(\mathbb{Z}/\langle p^k \rangle\) are \(u\)-\(S\)-torsion, hence \(u\)-\(S\)-absolutely pure. However, their direct product \(M = \prod_{k=1}^{\infty} \mathbb{Z}/\langle p^k \rangle\) is not \(u\)-\(S\)-injective, and therefore not \(u\)-\(S\)-absolutely pure by Theorem \ref{thm:usnoetherianiffapinjective} (since \(\mathbb{Z}\) is not \(u\)-\(S\)-Noetherian for this \(S\)). Similarly, their direct sum is also not \(u\)-\(S\)-absolutely pure.
	\end{example}
	
	A further example shows that the element \(s\) in the definition must be uniform; a module can have \(\Ext_R^1(N, E)\) \(u\)-\(S\)-torsion for all finitely presented \(N\) without being \(u\)-\(S\)-absolutely pure.
	
	\begin{example}\cite[Example 3.9]{Zhang2023d}
		\label{ex:extcyclefinebutnotap}
		Let \(R = \mathbb{Z}\), \(S = \{p^n \mid n \ge 0\}\), and \(J_p\) the \(p\)-adic integers. Then \(\Ext_R^1(N, J_p)\) is \(u\)-\(S\)-torsion for any finitely presented \(R\)-module \(N\). However, \(J_p\) is not \(u\)-\(S\)-injective (see Example \ref{ex:extcycleinsufficient}), and thus not \(u\)-\(S\)-absolutely pure by Theorem \ref{thm:usnoetherianiffapinjective} (since \(\mathbb{Z}\) is not \(u\)-\(S\)-Noetherian).
	\end{example}
	
	\section{Homological Dimensions: Uniformly S-Weak Global and Uniformly S-Global Dimensions}
	
	With the uniform analogues of flat, projective, and injective modules defined, we can now introduce the corresponding homological dimensions. These dimensions measure how far a ring is from being ``\(u\)-\(S\)-von Neumann regular'' or ``\(u\)-\(S\)-semisimple''.
	
	\subsection{Uniformly S-Flat Dimensions and Uniformly S-Weak Global Dimensions}
	
	The \(u\)-\(S\)-flat dimension of a module is defined as the length of the shortest \(u\)-\(S\)-flat resolution.
	
	\begin{definition}\cite[Definition 2.1]{Zhang2023b}
		\label{def:usflatdim}
		Let \(R\) be a ring, \(S\) a multiplicative subset of \(R\) and \(M\) an \(R\)-module. We write \(u\)-\(S\)-\(\operatorname{fd}_R(M) \le n\) if there exists a \(u\)-\(S\)-exact sequence of \(R\)-modules
		\[
		0 \to F_n \to \dots \to F_1 \to F_0 \to M \to 0
		\]
		where each \(F_i\) is \(u\)-\(S\)-flat for \(i = 0, \dots, n\). If such a finite \(u\)-\(S\)-flat \(u\)-\(S\)-resolution does not exist, then we say \(u\)-\(S\)-\(\operatorname{fd}_R(M) = \infty\); otherwise, define \(u\)-\(S\)-\(\operatorname{fd}_R(M) = n\) if \(n\) is the length of the shortest such resolution.
	\end{definition}
	
	Trivially, \(u\)-\(S\)-\(\operatorname{fd}_R(M) \le \operatorname{fd}_R(M)\). A module is \(u\)-\(S\)-flat iff its \(u\)-\(S\)-flat dimension is 0.
	
	The following lemma shows that the dimension is invariant under \(u\)-\(S\)-isomorphisms.
	
	\begin{lemma}\cite[Lemma 2.2]{Zhang2023b}
		\label{lem:usfdisoinv}
		Let \(R\) be a ring and \(S\) a multiplicative subset of \(R\). If \(A\) is \(u\)-\(S\)-isomorphic to \(B\), then \(u\)-\(S\)-\(\operatorname{fd}_R(A) = u\)-\(S\)-\(\operatorname{fd}_R(B)\).
	\end{lemma}
	
	The key characterization of \(u\)-\(S\)-flat dimension uses the Tor functor.
	
	\begin{proposition}\cite[Proposition 2.3]{Zhang2023b}
		\label{prop:usfdchar}
		Let \(R\) be a ring and \(S\) a multiplicative subset of \(R\). The following statements are equivalent for an \(R\)-module \(M\):
		\begin{enumerate}
			\item \(u\)-\(S\)-\(\operatorname{fd}_R(M) \le n\);
			\item \(\Tor_{n+k}^R(M,N)\) is \(u\)-\(S\)-torsion for all \(R\)-modules \(N\) and all \(k > 0\);
			\item \(\Tor_{n+1}^R(M,N)\) is \(u\)-\(S\)-torsion for all \(R\)-modules \(N\);
			\item there exists \(s \in S\) such that \(s\Tor_{n+1}^R(M,R/I) = 0\) for all ideals \(I\) of \(R\);
			\item if \(0 \to F_n \to \dots \to F_1 \to F_0 \to M \to 0\) is a \(u\)-\(S\)-exact sequence, where \(F_0, \dots, F_{n-1}\) are \(u\)-\(S\)-flat, then \(F_n\) is \(u\)-\(S\)-flat;
			\item[8.] if \(0 \to F_n \to \dots \to F_1 \to F_0 \to M \to 0\) is an exact sequence, where \(F_0, \dots, F_{n-1}\) are flat, then \(F_n\) is \(u\)-\(S\)-flat;
			\item[10.] there exists an exact sequence \(0 \to F_n \to \dots \to F_1 \to F_0 \to M \to 0\), where \(F_0, \dots, F_{n-1}\) are flat and \(F_n\) is \(u\)-\(S\)-flat.
		\end{enumerate}
	\end{proposition}
	
	A straightforward corollary relates dimensions for different multiplicative sets.
	
	\begin{corollary} \cite[Corollary 2.4]{Zhang2023b}
		\label{cor:usfdcomparison}
		Let \(R\) be a ring and \(S' \subseteq S\) multiplicative subsets of \(R\). Let \(M\) be an \(R\)-module. Then \(u\)-\(S\)-\(\operatorname{fd}_R(M) \le u\)-\(S'\)-\(\operatorname{fd}_R(M)\).
	\end{corollary}
	
	Localizing at the uniform element gives an upper bound on the classical flat dimension of the localized module.
	
	\begin{corollary}\cite[Corollary 2.5]{Zhang2023b}
		\label{cor:usfdloc}
		Let \(R\) be a ring, \(S\) a multiplicative subset of \(R\) and \(M\) an \(R\)-module. If \(u\)-\(S\)-\(\operatorname{fd}_R(M) \le n\), then there exists an element \(s \in S\) such that \(\operatorname{fd}_{R_s}(M_s) \le n\).
	\end{corollary}
	
	A more precise relationship holds when \(S\) is finite.
	
	\begin{corollary}\cite[Corollary 2.6]{Zhang2023b}
		\label{cor:usfdlocfinite}
		Let \(R\) be a ring and \(S\) a multiplicative subset of \(R\). Let \(M\) be an \(R\)-module. Then \(u\)-\(S\)-\(\operatorname{fd}_R(M) \ge \operatorname{fd}_{R_S} M_S\). Moreover, if \(S\) is composed of finite elements, then \(u\)-\(S\)-\(\operatorname{fd}_R(M) = \operatorname{fd}_{R_S} M_S\).
	\end{corollary}
	
	The \(u\)-\(S\)-flat dimension behaves like its classical counterpart in short exact sequences.
	
	\begin{proposition}\cite[Proposition 2.7]{Zhang2023b}
		\label{prop:usfdexact}
		Let \(R\) be a ring and \(S\) a multiplicative subset of \(R\). Let \(0 \to A \to B \to C \to 0\) be a \(u\)-\(S\)-exact sequence of \(R\)-modules. Then the following assertions hold.
		\begin{enumerate}
			\item \(u\)-\(S\)-\(\operatorname{fd}_R(C) \le 1 + \max\{u\)-\(S\)-\(\operatorname{fd}_R(A), u\)-\(S\)-\(\operatorname{fd}_R(B)\}\).
			\item If \(u\)-\(S\)-\(\operatorname{fd}_R(B) < u\)-\(S\)-\(\operatorname{fd}_R(C)\), then \(u\)-\(S\)-\(\operatorname{fd}_R(A) = u\)-\(S\)-\(\operatorname{fd}_R(C) - 1 > u\)-\(S\)-\(\operatorname{fd}_R(B)\).
		\end{enumerate}
	\end{proposition}
	
	A new local characterization of the classical flat dimension is obtained using the uniform dimensions at prime and maximal ideals.
	
	\begin{proposition}\cite[Proposition 2.8]{Zhang2023b}
		\label{prop:fdlocalus}
		Let \(R\) be a ring and \(M\) an \(R\)-module. Then
		\[
		\operatorname{fd}_R(M) = \sup \{u\mbox{-}\mathfrak{p}\mbox{-}\operatorname{fd}_R(M) \mid \mathfrak{p} \in \Spec(R)\} = \sup \{u\mbox{-}\mathfrak{m}\mbox{-}\operatorname{fd}_R(M) \mid \mathfrak{m} \in \Max(R)\}.
		\]
	\end{proposition}
	
	The \(u\)-\(S\)-weak global dimension of a ring is the supremum of the \(u\)-\(S\)-flat dimensions of all modules.
	
	\begin{definition}\cite[Definition 3.1]{Zhang2023b}
		\label{def:uswgd}
		The \textbf{\(u\)-\(S\)-weak global dimension} of a ring \(R\) is defined by \(u\)-\(S\)-\(\operatorname{w.gl.dim}(R) = \sup \{u\)-\(S\)-\(\operatorname{fd}_R(M) \mid M \text{ is an } R\text{-module} \}\).
	\end{definition}
	
	Obviously, \(u\)-\(S\)-\(\operatorname{w.gl.dim}(R) \le \operatorname{w.gl.dim}(R)\). The next result gives several equivalent characterizations.
	
	\begin{proposition}\cite[Proposition 3.2]{Zhang2023b}
		\label{prop:uswgdchar}
		Let \(R\) be a ring and \(S\) a multiplicative subset of \(R\). The following statements are equivalent for \(R\):
		\begin{enumerate}
			\item \(u\)-\(S\)-\(\operatorname{w.gl.dim}(R) \le n\);
			\item \(u\)-\(S\)-\(\operatorname{fd}_R(M) \le n\) for all \(R\)-modules \(M\);
			\item \(\Tor_{n+k}^R(M,N)\) is \(u\)-\(S\)-torsion for all \(R\)-modules \(M,N\) and all \(k > 0\);
			\item \(\Tor_{n+1}^R(M,N)\) is \(u\)-\(S\)-torsion for all \(R\)-modules \(M,N\);
			\item there exists an element \(s \in S\) such that \(s\Tor_{n+1}^R(R/I, R/J)\) for any ideals \(I\) and \(J\) of \(R\).
		\end{enumerate}
	\end{proposition}
	
	A direct corollary shows the relationship between dimensions for different multiplicative sets.
	
	\begin{corollary}\cite[Corollary 3.3]{Zhang2023b}
		\label{cor:uswgdcomparison}
		Let \(R\) be a ring and \(S' \subseteq S\) multiplicative subsets of \(R\). Then \(u\)-\(S\)-\(\operatorname{w.gl.dim}(R) \le u\)-\(S'\)-\(\operatorname{w.gl.dim}(R)\).
	\end{corollary}
	
	Localizing at the uniform element yields a bound on the classical weak global dimension.
	
	\begin{corollary}\cite[Corollary 3.4]{Zhang2023b}
		\label{cor:uswgdloc}
		Let \(R\) be a ring and \(S\) a multiplicative subset of \(R\). If \(u\)-\(S\)-\(\operatorname{w.gl.dim}(R) \le n\), then there exists an element \(s \in S\) such that \(\operatorname{w.gl.dim}(R_s) \le n\).
	\end{corollary}
	
	A more precise relationship holds when \(S\) is finite.
	
	\begin{corollary}\cite[Corollary 3.5]{Zhang2023b}
		\label{cor:uswgdlocfinite}
		Let \(R\) be a ring and \(S\) a multiplicative subset of \(R\). Then \(u\)-\(S\)-\(\operatorname{w.gl.dim}(R) \le \operatorname{w.gl.dim}(R_S)\). Moreover, if \(S\) is composed of finite elements, then \(u\)-\(S\)-\(\operatorname{w.gl.dim}(R) = \operatorname{w.gl.dim}(R_S)\).
	\end{corollary}
	
	An example shows that the inequality can be strict.
	
	\begin{example}\cite[Example 3.6]{Zhang2023b}
		\label{ex:uswgdneqwgdloc}
		Let \(R = k[x_1, \dots, x_{n+1}]\) be a polynomial ring over a field \(k\). Set \(S = k[x_1] \setminus \{0\}\). Then \(R_S = k(x_1)[x_2, \dots, x_{n+1}]\), and \(\operatorname{w.gl.dim}(R_S) = n\). However, for any \(s \in S\), \(\operatorname{w.gl.dim}(R_s) = n+1\). Therefore, \(u\)-\(S\)-\(\operatorname{w.gl.dim}(R) \ge n+1 > \operatorname{w.gl.dim}(R_S)\).
	\end{example}
	
	A local characterization of the classical weak global dimension follows.
	
	\begin{corollary} \cite[Corollary 3.7]{Zhang2023b}
		\label{cor:wgdlocalus}
		Let \(R\) be a ring. Then
		\[
		\operatorname{w.gl.dim}(R) = \sup \{u\mbox{-}\mathfrak{p}\mbox{-}\operatorname{w.gl.dim}(R) \mid \mathfrak{p} \in \Spec(R)\} = \sup \{u\mbox{-}\mathfrak{m}\mbox{-}\operatorname{w.gl.dim}(R) \mid \mathfrak{m} \in \Max(R)\}.
		\]
	\end{corollary}
	
	Rings with \(u\)-\(S\)-weak global dimension 0 are precisely the \(u\)-\(S\)-von Neumann regular rings, a result from \cite{Zhang2022}.
	
	\begin{corollary}\cite[Corollary 3.8]{Zhang2023b}
		\label{cor:uswgd0iffusvnr}
		Let \(R\) be a ring and \(S\) a multiplicative subset of \(R\). The following assertions are equivalent:
		\begin{enumerate}
			\item \(R\) is a \(u\)-\(S\)-von Neumann regular ring;
			\item for any \(R\)-module \(M\) and \(N\), there exists \(s \in S\) such that \(s\Tor_1^R(M,N) = 0\);
			\item there exists \(s \in S\) such that \(s\Tor_1^R(R/I, R/J) = 0\) for any ideals \(I\) and \(J\) of \(R\);
			\item any \(R\)-module is \(u\)-\(S\)-flat;
			\item \(u\)-\(S\)-\(\operatorname{w.gl.dim}(R) = 0\).
		\end{enumerate}
	\end{corollary}
	
	Rings with \(u\)-\(S\)-weak global dimension at most 1 are characterized by the property that submodules of \(u\)-\(S\)-flat modules are \(u\)-\(S\)-flat.
	
	\begin{proposition}\cite[Proposition 3.9]{Zhang2023b}
		\label{prop:uswgd1char}
		Let \(R\) be a ring and \(S\) a multiplicative subset of \(R\). The following assertions are equivalent:
		\begin{enumerate}
			\item \(u\)-\(S\)-\(\operatorname{w.gl.dim}(R) \le 1\);
			\item any submodule of \(u\)-\(S\)-flat modules is \(u\)-\(S\)-flat;
			\item any submodule of flat modules is \(u\)-\(S\)-flat;
			\item \(\Tor_2^R(M,N)\) is \(u\)-\(S\)-torsion for all \(R\)-modules \(M,N\);
			\item there exists an element \(s \in S\) such that \(s\Tor_2^R(R/I, R/J) = 0\) for any ideals \(I,J\) of \(R\).
		\end{enumerate}
	\end{proposition}

	The following example demonstrates that a ring can have finite \(u\)-\(S\)-weak global dimension while having infinite classical weak global dimension.
	
	\begin{example}\cite[Example 3.11]{Zhang2023b}
		\label{ex:uswgdfinitewgdinfinite}
		Let \(A\) be a ring with \(\operatorname{w.gl.dim}(A) = 1\), \(T = A \times A\), and \(s = (1,0) \in T\). Let \(R = T[x]/\langle sx, x^2 \rangle\) and \(S = \{1, s\}\). Then \(u\)-\(S\)-\(\operatorname{w.gl.dim}(R) = 1\) but \(\operatorname{w.gl.dim}(R) = \infty\).
	\end{example}
	
	\subsection{Uniformly S-Projective (Injective) Dimensions, and Uniformly S-Global Dimensions}
	
	The uniform analogues of projective and injective dimensions are defined similarly.
	
	\begin{definition}\cite[Definition 3.1]{Zhang2023c}
		\label{def:uspddim}
		Let \(R\) be a ring, \(S\) a multiplicative subset of \(R\) and \(M\) an \(R\)-module. We write \(u\)-\(S\)-\(\operatorname{pd}_R(M) \le n\) if there exists a \(u\)-\(S\)-exact sequence of \(R\)-modules
		\[
		0 \to F_n \to \dots \to F_1 \to F_0 \to M \to 0
		\]
		where each \(F_i\) is \(u\)-\(S\)-projective for \(i = 0, \dots, n\). If such a finite \(u\)-\(S\)-projective \(u\)-\(S\)-resolution does not exist, then we say \(u\)-\(S\)-\(\operatorname{pd}_R(M) = \infty\); otherwise, define \(u\)-\(S\)-\(\operatorname{pd}_R(M) = n\) if \(n\) is the length of the shortest such resolution. The \(u\)-\(S\)-injective dimension \(u\)-\(S\)-\(\operatorname{id}_R(M)\) is defined dually.
	\end{definition}
	
	These dimensions are invariant under \(u\)-\(S\)-isomorphisms.
	
	\begin{lemma}\cite[Lemma 3.2]{Zhang2023c}
		\label{lem:usdimisoinv}
		Let \(R\) be a ring, \(S\) a multiplicative subset of \(R\). If \(A\) is \(u\)-\(S\)-isomorphic to \(B\), then \(u\)-\(S\)-\(\operatorname{pd}_R(A) = u\)-\(S\)-\(\operatorname{pd}_R(B)\) and \(u\)-\(S\)-\(\operatorname{id}_R(A) = u\)-\(S\)-\(\operatorname{id}_R(B)\).
	\end{lemma}
	
	The characterization of \(u\)-\(S\)-projective dimension uses the Ext functor.
	
	\begin{proposition}\cite[Proposition 3.3]{Zhang2023c}
		\label{prop:uspddimchar}
		Let \(R\) be a ring and \(S\) a multiplicative subset of \(R\). The following statements are equivalent for an \(R\)-module \(M\):
		\begin{enumerate}
			\item \(u\)-\(S\)-\(\operatorname{pd}_R(M) \le n\);
			\item \(\Ext_{R}^{n+k}(M,N)\) is \(u\)-\(S\)-torsion for all \(R\)-modules \(N\) and all \(k > 0\);
			\item \(\Ext_{R}^{n+1}(M,N)\) is \(u\)-\(S\)-torsion for all \(R\)-modules \(N\);
			\item if \(0 \to F_n \to \dots \to F_1 \to F_0 \to M \to 0\) is a \(u\)-\(S\)-exact sequence, where \(F_0, \dots, F_{n-1}\) are \(u\)-\(S\)-projective, then \(F_n\) is \(u\)-\(S\)-projective;
			\item[7.] if \(0 \to F_n \to \dots \to F_1 \to F_0 \to M \to 0\) is an exact sequence, where \(F_0, \dots, F_{n-1}\) are projective, then \(F_n\) is \(u\)-\(S\)-projective;
			\item[9.] there exists an exact sequence \(0 \to F_n \to \dots \to F_1 \to F_0 \to M \to 0\), where \(F_0, \dots, F_{n-1}\) are projective and \(F_n\) is \(u\)-\(S\)-projective.
		\end{enumerate}
	\end{proposition}
	
	A similar characterization holds for the \(u\)-\(S\)-injective dimension.
	
	\begin{proposition}\cite[Proposition 3.4]{Zhang2023c}
		\label{prop:usiddimchar}
		Let \(R\) be a ring and \(S\) a multiplicative subset of \(R\). The following statements are equivalent for an \(R\)-module \(M\):
		\begin{enumerate}
			\item \(u\)-\(S\)-\(\operatorname{id}_R(M) \le n\);
			\item \(\Ext_{R}^{n+k}(N,M)\) is \(u\)-\(S\)-torsion for all \(R\)-modules \(N\) and all \(k > 0\);
			\item \(\Ext_{R}^{n+1}(N,M)\) is \(u\)-\(S\)-torsion for all \(R\)-modules \(N\);
			\item if \(0 \to M \to E_0 \to \dots \to E_{n-1} \to E_n \to 0\) is a \(u\)-\(S\)-exact sequence, where \(E_0, \dots, E_{n-1}\) are \(u\)-\(S\)-injective, then \(E_n\) is \(u\)-\(S\)-injective;
			\item[7.] if \(0 \to M \to E_0 \to \dots \to E_{n-1} \to E_n \to 0\) is an exact sequence, where \(E_0, \dots, E_{n-1}\) are injective, then \(E_n\) is \(u\)-\(S\)-injective;
			\item[9.] there exists an exact sequence \(0 \to M \to E_0 \to \dots \to E_{n-1} \to E_n \to 0\), where \(E_0, \dots, E_{n-1}\) are injective and \(E_n\) is \(u\)-\(S\)-injective.
		\end{enumerate}
	\end{proposition}
	
	For \(u\)-\(S\)-Noetherian rings, the \(u\)-\(S\)-projective and \(u\)-\(S\)-flat dimensions of \(S\)-finite modules coincide.
	
	\begin{proposition}\cite[Proposition 3.6]{Zhang2023c}
		\label{prop:uspdequalsfdforusnoeth}
		Let \(R\) be a ring and \(S\) a multiplicative subset of \(R\). If \(R\) is a \(u\)-\(S\)-Noetherian ring, then the following statements hold.
		\begin{enumerate}
			\item If \(M\) is an \(S\)-finite \(R\)-module, then there is a \(u\)-\(S\)-exact sequence
			\[
			\dots \to F_n \to \dots \to F_1 \to F_0 \to M \to 0
			\]
			with each \(F_n\) \(S\)-finite \(u\)-\(S\)-projective.
			\item If \(M\) is an \(S\)-finite \(R\)-module, then \(u\)-\(S\)-\(\operatorname{pd}_R(M) = u\)-\(S\)-\(\operatorname{fd}_R(M)\).
		\end{enumerate}
	\end{proposition}

	The behavior of these dimensions in \(u\)-\(S\)-exact sequences is standard.
	
	\begin{proposition}\cite[Proposition 3.10]{Zhang2023c}
		\label{prop:usdimexact}
		Let \(R\) be a ring and \(S\) a multiplicative subset of \(R\). Let \(0 \to A \to B \to C \to 0\) be a \(u\)-\(S\)-exact sequence of \(R\)-modules. Then the following statements hold.
		\begin{enumerate}
			\item \(u\)-\(S\)-\(\operatorname{pd}_R(C) \le 1 + \max\{u\)-\(S\)-\(\operatorname{pd}_R(A), u\)-\(S\)-\(\operatorname{pd}_R(B)\}\).
			\item If \(u\)-\(S\)-\(\operatorname{pd}_R(B) < u\)-\(S\)-\(\operatorname{pd}_R(C)\), then \(u\)-\(S\)-\(\operatorname{pd}_R(A) = u\)-\(S\)-\(\operatorname{pd}_R(C) - 1 > u\)-\(S\)-\(\operatorname{pd}_R(B)\).
			\item \(u\)-\(S\)-\(\operatorname{id}_R(A) \le 1 + \max\{u\)-\(S\)-\(\operatorname{id}_R(B), u\)-\(S\)-\(\operatorname{id}_R(C)\}\).
			\item If \(u\)-\(S\)-\(\operatorname{id}_R(B) < u\)-\(S\)-\(\operatorname{id}_R(A)\), then \(u\)-\(S\)-\(\operatorname{id}_R(C) = u\)-\(S\)-\(\operatorname{id}_R(A) - 1 > u\)-\(S\)-\(\operatorname{id}_R(B)\).
		\end{enumerate}
	\end{proposition}
	
	For \(u\)-\(S\)-split sequences, the dimension of the middle term is the maximum of the dimensions of the ends.
	
	\begin{proposition}\cite[Proposition 3.11]{Zhang2023c}
		\label{prop:usdimsplit}
		Let \(0 \to A \to B \to C \to 0\) be a \(u\)-\(S\)-split \(u\)-\(S\)-exact sequence of \(R\)-modules. Then the following statements hold.
		\begin{enumerate}
			\item \(u\)-\(S\)-\(\operatorname{pd}_R(B) = \max\{u\)-\(S\)-\(\operatorname{pd}_R(A), u\)-\(S\)-\(\operatorname{pd}_R(C)\}\).
			\item \(u\)-\(S\)-\(\operatorname{id}_R(B) = \max\{u\)-\(S\)-\(\operatorname{id}_R(A), u\)-\(S\)-\(\operatorname{id}_R(C)\}\).
		\end{enumerate}
	\end{proposition}
	
	A new local characterization of the classical projective and injective dimensions is obtained.
	
	\begin{proposition}\cite[Proposition 3.12]{Zhang2023c}
		\label{prop:pdidlocalus}
		Let \(R\) be a ring and \(M\) an \(R\)-module. Then
		\[
		\operatorname{pd}_R(M) = \sup \{u\mbox{-}\mathfrak{p}\mbox{-}\operatorname{pd}_R(M) \mid \mathfrak{p}\in \Spec(R)\} = \sup \{u\mbox{-}\mathfrak{m}\mbox{-}\operatorname{pd}_R(M)\mid \mathfrak{m} \in \Max(R)\},
		\]
		and
		\[
		\operatorname{id}_R(M) = \sup \{u\mbox{-}\mathfrak{p}\mbox{-}\operatorname{id}_R(M) \mid \mathfrak{p} \in \Spec(R)\} = \sup \{u\mbox{-}\mathfrak{m}\mbox{-}\operatorname{id}_R(M) \mid \mathfrak{m} \in \Max(R)\}.
		\]
	\end{proposition}
	
	The \(u\)-\(S\)-global dimension of a ring is the supremum of the \(u\)-\(S\)-projective dimensions of all modules, which also equals the supremum of the \(u\)-\(S\)-injective dimensions.
	
	\begin{definition}\cite[Definition 3.13]{Zhang2023c}
		\label{def:usgldim}
		The \textbf{\(u\)-\(S\)-global dimension} of a ring \(R\) is defined by
		\[
		u\mbox{-}S\mbox{-}\operatorname{gl.dim}(R) = \sup \{u\mbox{-}S\mbox{-}\operatorname{pd}_R(M) \mid M \text{ is an } R\text{-module} \}.
		\]
	\end{definition}
	
	The following theorem provides several equivalent characterizations, mirroring those for the classical global dimension.
	
	\begin{proposition}\cite[Proposition 3.14]{Zhang2023c}
		\label{prop:usgldimchar}
		Let \(R\) be a ring and \(S\) a multiplicative subset of \(R\). The following statements are equivalent for \(R\):
		\begin{enumerate}
			\item \(u\)-\(S\)-\(\operatorname{gl.dim}(R) \le n\);
			\item \(u\)-\(S\)-\(\operatorname{pd}_R(M) \le n\) for all \(R\)-modules \(M\);
			\item \(\Ext_{R}^{n+k}(M,N)\) is \(u\)-\(S\)-torsion for all \(R\)-modules \(M,N\) and all \(k > 0\);
			\item \(\Ext_{R}^{n+1}(M,N)\) is \(u\)-\(S\)-torsion for all \(R\)-modules \(M,N\);
			\item \(u\)-\(S\)-\(\operatorname{id}_R(M) \le n\) for all \(R\)-modules \(M\).
		\end{enumerate}
	\end{proposition}
	
	A direct corollary shows that the \(u\)-\(S\)-global dimension bounds the \(u\)-\(S\)-weak global dimension from above.
	
	\begin{corollary}\cite[Corollary 3.15]{Zhang2023c}
		\label{cor:usgldimgequswgd}
		Let \(R\) be a ring, \(S\) a multiplicative subset of \(R\). Then \(u\)-\(S\)-\(\operatorname{gl.dim}(R) \ge u\)-\(S\)-\(\operatorname{w.gl.dim}(R)\).
	\end{corollary}
	
	The dimension behaves predictably under changes of the multiplicative set.
	
	\begin{corollary} \cite[Corollary 3.16]{Zhang2023c}
		\label{cor:usgldimcomparison}
		Let \(R\) be a ring, \(S' \subseteq S\) be multiplicative subsets of \(R\) and \(S^{*}\) be the saturation of \(S\). Then \(u\)-\(S\)-\(\operatorname{gl.dim}(R) \le u\)-\(S'\)-\(\operatorname{gl.dim}(R)\) and \(u\)-\(S\)-\(\operatorname{gl.dim}(R) = u\)-\(S^{*}\)-\(\operatorname{gl.dim}(R)\).
	\end{corollary}
	
	For a finite direct product of rings, the \(u\)-\(S\)-global dimension is the supremum of the dimensions of the factors.
	
	\begin{corollary}\cite[Corollary 3.17]{Zhang2023c}
		\label{cor:usgldimproduct}
		Let \(R_i\) be a ring and \(S_i\) be a multiplicative subset of \(R_i\) \((i = 1, \dots, n)\). Set \(R = R_1 \times \dots \times R_n\) and \(S = S_1 \times \dots \times S_n\) a multiplicative subset of \(R_i\). Then \(u\)-\(S\)-\(\operatorname{gl.dim}(R) = \sup_{1 \le i \le n}\{u\)-\(S_i\)-\(\operatorname{gl.dim}(R_i)\}\).
	\end{corollary}
	
	An example shows that the classical global dimension and the \(u\)-\(S\)-global dimension can be quite different.
	
	\begin{example}\cite[Example 3.18]{Zhang2023c}
		\label{ex:usgldimneqgldim}
		Let \(R_1\) be a ring with \(\operatorname{gl.dim}(R_1) = n\) and \(R_2\) be a ring with \(\operatorname{gl.dim}(R_2) = m\). Set \(R = R_1 \times R_2\) and \(S = \{(1,1),(1,0)\}\). Then \(\operatorname{gl.dim}(R) = \max\{m,n\}\), but \(u\)-\(S\)-\(\operatorname{gl.dim}(R) = n\).
	\end{example}
	
	A local characterization of the classical global dimension follows.
	
	\begin{corollary}\cite[Corollary 3.19]{Zhang2023c}
		\label{cor:gldimlocalus}
		Let \(R\) be a ring. Then
		\[
		\operatorname{gl.dim}(R) = \sup \{u\mbox{-}\mathfrak{p}\mbox{-}\operatorname{gl.dim}(R) \mid \mathfrak{p} \in \Spec(R)\} = \sup \{u\mbox{-}\mathfrak{m}\mbox{-}\operatorname{gl.dim}(R) \mid \mathfrak{m} \in \Max(R)\}.
		\]
	\end{corollary}
	
	Rings with \(u\)-\(S\)-global dimension 0 are precisely the \(u\)-\(S\)-semisimple rings, a concept introduced in \cite{Zhang2023a}.
	
	\begin{corollary} \cite[Corollary 3.20]{Zhang2023c}
		\label{cor:usgldim0iffussemisimple}
		Let \(R\) be a ring and \(S\) a multiplicative subset of \(R\). The following statements are equivalent:
		\begin{enumerate}
			\item \(R\) is a \(u\)-\(S\)-semisimple ring;
			\item every \(R\)-module is \(u\)-\(S\)-semisimple;
			\item every \(R\)-module is \(u\)-\(S\)-projective;
			\item every \(R\)-module is \(u\)-\(S\)-injective;
			\item \(u\)-\(S\)-\(\operatorname{gl.dim}(R) = 0\).
		\end{enumerate}
	\end{corollary}
	
	\subsection{Change of Rings Theorems for Uniformly S-Global Dimensions}
	
	Understanding how homological dimensions behave under ring homomorphisms and extensions is a central theme in homological algebra. The uniform \(S\)-theory is no exception. The following proposition provides an upper bound for the \(u\)-\(S\)-projective dimension of a module when changing rings.
	
	\begin{proposition}\cite[Proposition 4.1]{Zhang2023c}
		\label{prop:uspdchangeofrings}
		Let \(\theta: R \to T\) be a ring homomorphism, \(S\) a multiplicative subset of \(R\). Suppose \(M\) is a \(T\)-module. Then
		\[
		u\mbox{-}S\mbox{-}\operatorname{pd}_R(M) \le u\mbox{-}\theta(S)\mbox{-}\operatorname{pd}_T(M) + u\mbox{-}S\mbox{-}\operatorname{pd}_R(T).
		\]
	\end{proposition}
	
	A special case is when we pass to a quotient ring by a non-zero-divisor. This yields a formula relating dimensions over the original ring and the quotient.
	
	\begin{proposition}\cite[Proposition 4.2]{Zhang2023c}
		\label{prop:uspdquotient}
		Let \(R\) be a ring, \(S\) a multiplicative subset of \(R\). Let \(a\) be a nonzero-divisor in \(R\) which does not divide any element in \(S\). Write \(\overline{R} = R/aR\) and \(\overline{S} = \{s + aR \in \overline{R} \mid s \in S\}\). Then the following statements hold.
		\begin{enumerate}
			\item Let \(M\) be a nonzero \(\overline{R}\)-module. If \(u\)\mbox{-}\(\overline{S}\)\mbox{-}\(\operatorname{pd}_{\overline{R}}(M) < \infty\), then
			\[
			u\mbox{-}S\mbox{-}\operatorname{pd}_R(M) = u\mbox{-}\overline{S}\mbox{-}\operatorname{pd}_{\overline{R}}(M) + 1.
			\]
			\item If \(u\)\mbox{-}\(\overline{S}\)\mbox{-}\(\operatorname{gl.dim}(\overline{R}) < \infty\), then
			\[
			u\mbox{-}S\mbox{-}\operatorname{gl.dim}(R) \ge u\mbox{-}\overline{S}\mbox{-}\operatorname{gl.dim}(\overline{R}) + 1.
			\]
		\end{enumerate}
	\end{proposition}
	
	A key result is the behavior of the \(u\)-\(S\)-projective dimension under polynomial extensions.
	
	\begin{lemma}\cite[Lemma 4.3]{Zhang2023c}
		\label{lem:usprojpolynomial}
		Let \(R\) be a ring, \(S\) a multiplicative subset of \(R\). Suppose \(T\) is an \(R\)-module and \(F\) is an \(R[x]\)-module. If \(P\) is \(u\)-\(S\)-projective over \(R[x]\), then \(P\) is \(u\)-\(S\)-projective over \(R\).
	\end{lemma}
	
	For the module \(M[x] = M \otimes_R R[x]\), the dimension is unchanged.
	
	\begin{proposition}\cite[Proposition 4.4]{Zhang2023c}
		\label{prop:uspdmxinvar}
		Let \(R\) be a ring, \(S\) a multiplicative subset of \(R\) and \(M\) an \(R\)-module. Then \(u\)-\(S\)-\(\operatorname{pd}_{R[x]}(M[x]) = u\)-\(S\)-\(\operatorname{pd}_R(M)\).
	\end{proposition}
	
We can prove the \(u\)-\(S\)-analogue of the Hilbert syzygy theorem. The result states that, under a mild condition, the \(u\)-\(S\)-global dimension increases by exactly one when passing to a polynomial ring.
	
	\begin{theorem}\cite[Theorem 4.6]{Zhang2023c}
		\label{thm:usgldimPolynomial}
		Let \(R\) be a ring and \(S\) a multiplicative subset of \(R\). Then
		\[
		u\mbox{-}S\mbox{-}\operatorname{gl.dim}(R[x]) = u\mbox{-}S\mbox{-}\operatorname{gl.dim}(R) + 1,
		\]
		provided that \(0 \notin S\).
	\end{theorem}
	
	A remark notes that the condition \(u\)-\(S\)-\(\operatorname{pd}_{R[x]}(R)=1\) is needed for the proof, and it holds in important cases like when \(S\) consists of units or for certain product rings. The theorem can fail if \(0 \in S\).
	
	\section{Uniformly S-version of Classical Rings}
	
	Building on the uniform homological algebra, several important classes of rings have been defined and characterized by uniform properties. These include \(u\)-\(S\)-von Neumann regular rings, \(u\)-\(S\)-semisimple rings, and \(u\)-\(S\)-Artinian rings.
	
	\subsection{Uniformly S-von Neumann Regular Rings}
	
	Von Neumann regular rings are characterized by the property that for every element \(a\), there exists an element \(r\) such that \(a = ra^2\). The uniform \(S\)-analogue requires a single element \(s \in S\) to work for all \(a\).
	
	\begin{definition}\cite[Definition 3.12]{Zhang2022}
		\label{def:usvnr}
		Let \(R\) be a ring and \(S\) a multiplicative subset of \(R\). \(R\) is called a \textbf{\(u\)-\(S\)-von Neumann regular ring} (abbreviates uniformly \(S\)-von Neumann regular ring) provided there exists an element \(s \in S\) satisfying that for any \(a \in R\) there exists \(r \in R\) such that \(sa = ra^2\).
	\end{definition}
	
	The following theorem provides a plethora of equivalent characterizations, linking this notion to \(u\)-\(S\)-flatness, torsion properties of \(\Tor\), and generation of ideals by idempotents.
	
	\begin{theorem}\cite[Theorem 3.13]{Zhang2022}
		\label{thm:usvnrchar}
		Let \(R\) be a ring and \(S\) a multiplicative subset of \(R\). The following statements are equivalent:
		\begin{enumerate}
			\item \(R\) is a \(u\)-\(S\)-von Neumann regular ring;
			\item for any \(R\)-module \(M\) and \(N\), there exists \(s \in S\) such that \(s\Tor_1^R(M,N) = 0\);
			\item there exists \(s \in S\) such that \(s\Tor_1^R(R/I, R/J) = 0\) for any ideals \(I\) and \(J\) of \(R\);
			\item there exists \(s \in S\) such that \(s\Tor_1^R(R/I, R/J) = 0\) for any \(S\)-finite ideals \(I\) and \(J\) of \(R\);
			\item there exists \(s \in S\) such that \(s\Tor_1^R(R/\langle a\rangle, R/\langle a\rangle) = 0\) for any element \(a \in R\);
			\item any \(R\)-module is \(u\)-\(S\)-flat;
			\item the class of all principal ideals of \(R\) is \(u\)-\(S\)-generated by idempotents;
			\item the class of all finitely generated ideals of \(R\) is \(u\)-\(S\)-generated by idempotents.
		\end{enumerate}
	\end{theorem}
	
	Localizing a \(u\)-\(S\)-von Neumann regular ring at \(S\) yields a classical von Neumann regular ring.
	
	\begin{corollary} \cite[Corollary 3.14]{Zhang2022}
		\label{cor:usvnrloc}
		Let \(R\) be a ring and \(S\) a multiplicative subset of \(R\). If \(R\) is a \(u\)-\(S\)-von Neumann regular ring, then \(R_S\) is a von Neumann regular ring. Consequently, any \(u\)-\(S\)-von Neumann regular ring is \(c\)-\(S\)-coherent.
	\end{corollary}
	
	However, a ring for which \(R_S\) is von Neumann regular need not be \(u\)-\(S\)-von Neumann regular.
	
	\begin{example}\cite[Example 3.15]{Zhang2022}
		\label{ex:locvnrnotusvnr}
		Let \(R = \mathbb{Z}\) and \(S = \mathbb{Z} \setminus \{0\}\). Then \(R_S = \mathbb{Q}\) is von Neumann regular. However, \(R\) is not \(u\)-\(S\)-von Neumann regular because \(\Tor_1^{\mathbb{Z}}(\mathbb{Q}/\mathbb{Z}, \mathbb{Z}_{(p)}/\mathbb{Z}) \cong \mathbb{Z}_{(p)}/\mathbb{Z}\) is not \(u\)-\(S\)-torsion.
	\end{example}
	
	When \(S\) is finite, the two properties coincide.
	
	\begin{corollary}\cite[Corollary 3.16]{Zhang2022}
		\label{cor:usvnrlocfinite}
		Let \(R\) be a ring. Let \(S\) be a multiplicative subset of \(R\) consisting of finite elements. Then \(R\) is a \(u\)-\(S\)-von Neumann regular ring if and only if \(R_S\) is a von Neumann regular ring.
	\end{corollary}
	
	If \(S\) consists of non-zero-divisors, then a \(u\)-\(S\)-von Neumann regular ring is actually von Neumann regular.
	
	\begin{proposition}\cite[Proposition 3.17]{Zhang2022}
		\label{prop:usvnrregularimpliesvnr}
		Let \(R\) be a ring and \(S\) a regular multiplicative subset of \(R\). Then \(R\) is \(u\)-\(S\)-von Neumann regular if and only if \(R\) is von Neumann regular.
	\end{proposition}
	
	A non-trivial example of a \(u\)-\(S\)-von Neumann regular ring that is not von Neumann regular shows that the regularity condition on \(S\) is essential.
	
	\begin{example} \cite[Example 3.18]{Zhang2022}
		\label{ex:usvnrnotvnr}
		Let \(T = \mathbb{Z}_2 \times \mathbb{Z}_2\) and \(s = (1,0) \in T\). Let \(R = T[x]/\langle sx, x^2 \rangle\) with \(x\) an indeterminate and \(S = \{1, \bar{s}\}\). Then \(R\) is a \(u\)-\(S\)-von Neumann regular ring, but \(R\) is not von Neumann regular (as it is not reduced).
	\end{example}
	
	A new local characterization of classical von Neumann regular rings is given by the uniform property at prime ideals.
	
	\begin{proposition}\cite[Proposition 3.19]{Zhang2022}
		\label{prop:vnrlocalus}
		Let \(R\) be a ring. Then the following statements are equivalent:
		\begin{enumerate}
			\item \(R\) is a von Neumann regular ring;
			\item \(R\) is a \(u\)-\(\mathfrak{p}\)-von Neumann regular ring for any \(\mathfrak{p} \in \Spec(R)\);
			\item \(R\) is a \(u\)-\(\mathfrak{m}\)-von Neumann regular ring for any \(\mathfrak{m} \in \Max(R)\).
		\end{enumerate}
	\end{proposition}
	
	\subsection{Uniformly S-Semisimple Rings}
	
	Semisimple rings are those where every module is projective. The uniform \(S\)-analogue relaxes this condition.
	
	\begin{definition} \cite[Definition 3.1]{Zhang2023a}
		\label{def:ussemisimplemodule}
		Let \(R\) be a ring and \(S\) a multiplicative subset of \(R\). An \(R\)-module \(M\) is called \textbf{\(u\)-\(S\)-semisimple} provided that any \(u\)-\(S\)-short exact sequence \(0 \to A \to M \to C \to 0\) is \(u\)-\(S\)-split.
	\end{definition}
	
	A ring is \(u\)-\(S\)-semisimple if all its free modules have this property.
	
	\begin{definition}\cite[Definition 3.4]{Zhang2023a}
		\label{def:ussemisimplering}
		Let \(R\) be a ring and \(S\) a multiplicative subset of \(R\). \(R\) is called a \textbf{\(u\)-\(S\)-semisimple ring} provided that any free \(R\)-module is \(u\)-\(S\)-semisimple.
	\end{definition}
	
	The following theorem provides a wealth of characterizations, showing that the property is incredibly strong.
	
	\begin{theorem} \cite[Theorem 3.5]{Zhang2023a}
		\label{thm:ussemisimplechar}
		Let \(R\) be a ring and \(S\) a multiplicative subset of \(R\). Then the following statements are equivalent:
		\begin{enumerate}
			\item \(R\) is a \(u\)-\(S\)-semisimple ring;
			\item any \(R\)-module is \(u\)-\(S\)-semisimple;
			\item any \(u\)-\(S\)-short exact sequence is \(u\)-\(S\)-split;
			\item any short exact sequence is \(u\)-\(S\)-split;
			\item \(\Ext_R^1(M,N)\) is \(u\)-\(S\)-torsion for any \(R\)-modules \(M\) and \(N\);
			\item any \(R\)-module is \(u\)-\(S\)-projective;
			\item any \(R\)-module is \(u\)-\(S\)-injective.
		\end{enumerate}
	\end{theorem}
	
	A \(u\)-\(S\)-semisimple ring is both \(u\)-\(S\)-Noetherian and \(u\)-\(S\)-von Neumann regular.
	
	\begin{corollary} \cite[Corollary 3.6]{Zhang2023a}
		\label{cor:ussemisimpleimpliesusnoethusvnr}
		Let \(R\) be a ring and \(S\) a multiplicative subset of \(R\). Suppose \(R\) is a \(u\)-\(S\)-semisimple ring. Then \(R\) is both \(u\)-\(S\)-Noetherian and \(u\)-\(S\)-von Neumann regular. Consequently, there exists an element \(s \in S\) such that for any ideal \(I\) of \(R\) there is an \(R\)-homomorphism \(f_I: R \to I\) satisfying \(f_I(i) = si\) for any \(i \in I\).
	\end{corollary}
	
	A ring can be \(u\)-\(S\)-semisimple as a module without being a \(u\)-\(S\)-semisimple ring.
	
	\begin{example}\cite[Example 3.7]{Zhang2023a}
		\label{ex:usmodssnotusringss}
		Let \(R = \mathbb{Z}\) and \(S = \mathbb{Z}\setminus\{0\}\). Then \(R\) is a \(u\)-\(S\)-semisimple \(\mathbb{Z}\)-module (for any ideal \(I\), pick \(s \in I\) and define \(f:R \to I\) by \(f(1)=s\)). However, \(R\) is not a \(u\)-\(S\)-semisimple ring because it is not \(u\)-\(S\)-von Neumann regular (see Example \ref{ex:locvnrnotusvnr}).
	\end{example}
	
	When \(S\) consists of non-zero-divisors, a \(u\)-\(S\)-semisimple ring is classically semisimple.
	
	\begin{proposition} \cite[Proposition 3.8]{Zhang2023a}
		\label{prop:ussemisimpleregularimpliesss}
		Let \(R\) be a ring and \(S\) a regular multiplicative subset of \(R\). Then \(R\) is a \(u\)-\(S\)-semisimple ring if and only if \(R\) is a semisimple ring.
	\end{proposition}
	
	The property is well-behaved under finite direct products.
	
	\begin{proposition} \cite[Proposition 3.10]{Zhang2023a}
		\label{prop:ussemisimpleproduct}
		Let \(R = R_1 \times R_2\) be direct product of rings \(R_1\) and \(R_2\) and \(S = S_1 \times S_2\) a direct product of multiplicative subsets. Then \(R\) is a \(u\)-\(S\)-semisimple ring if and only if \(R_i\) is a \(u\)-\(S_i\)-semisimple ring for each \(i = 1,2\).
	\end{proposition}
	
	An example shows that a ring can be \(u\)-\(S\)-semisimple without being semisimple, by taking a product with a non-semisimple ring and a suitable \(S\).
	
	\begin{example}\cite[Example 3.11]{Zhang2023a}
		\label{ex:ussemisimplenotss}
		Let \(R_1\) be a semisimple ring and \(R_2\) a non-semisimple ring. Let \(R = R_1 \times R_2\) and \(S = \{(1,1), (1,0)\}\). Then \(R\) is \(u\)-\(S\)-semisimple but not semisimple.
	\end{example}
	
	A local characterization of classical semisimple rings follows.
	
	\begin{proposition} \cite[Proposition 3.12]{Zhang2023a}
		\label{prop:sslocalus}
		Let \(R\) be a ring. Then the following statements are equivalent:
		\begin{enumerate}
			\item \(R\) is a semisimple ring;
			\item \(R\) is a \(u\)-\(\mathfrak{p}\)-semisimple ring for any \(\mathfrak{p} \in \Spec(R)\);
			\item \(R\) is a \(u\)-\(\mathfrak{m}\)-semisimple ring for any \(\mathfrak{m} \in \Max(R)\).
		\end{enumerate}
	\end{proposition}

	\subsection{Uniformly S-Multiplication Modules and Rings}
	
	Multiplication modules are those where every submodule can be written as a product of an ideal with the module. The \(S\)-version allows for a factor \(s\) that depends on the submodule. The uniform version fixes this \(s\).
	
	\begin{definition}\cite[Definition 1]{Qi2025}
		\label{def:usmultmodule}
		Let \(M\) be an \(R\)-module and let \(S\) be a multiplicative subset of \(R\). Then \(M\) is called a \textbf{\(u\)-\(S\)-multiplication} module (with respect to \(s\)) if there exists an element \(s \in S\) such that, for each submodule \(N\) of \(M\), there is an ideal \(I\) of \(R\) satisfying \(sN \subseteq IM \subseteq N\).
	\end{definition}
	
	Equivalently, \(M\) is a \(u\)-\(S\)-multiplication module iff there exists \(s \in S\) such that for every submodule \(N\), \(sN \subseteq (N:_R M)M \subseteq N\).
	
	The property is well-behaved under finite direct products.
	
	\begin{proposition} \cite[Proposition 1]{Qi2025}
		\label{prop:usmultmodproduct}
		Let \(M_i\) be an \(R_i\)-module and let \(S_i \subseteq R_i\) be a multiplicative subset \((i = 1,2)\). Set \(R = R_1 \times R_2, S = S_1 \times S_2\), and \(M = M_1 \times M_2\). Then \(M\) is a \(u\)-\(S\)-multiplication module if and only if \(M_1\) is a \(u\)-\(S_1\)-multiplication module and \(M_2\) is a \(u\)-\(S_2\)-multiplication module.
	\end{proposition}
	
	An example shows that an \(S\)-multiplication module need not be a \(u\)-\(S\)-multiplication module.
	
	\begin{example}\cite[Example 1]{Qi2025}
		\label{ex:smultnotusmult}
		Consider the \(\mathbb{Z}\)-module \(E(p) = \{\frac{r}{p^m} + \mathbb{Q} \in \mathbb{Q}/\mathbb{Z}\}\), with \(S = \{p^n : n \in \mathbb{N} \cup \{0\}\}\). Then \(E(p)\) is an \(S\)-multiplication module but not a \(u\)-\(S\)-multiplication module.
	\end{example}
	
	The property is invariant under taking the saturation of \(S\) and under \(u\)-\(S\)-isomorphisms.
	
	\begin{proposition}\cite[Proposition 2]{Qi2025}
		\label{prop:usmultmodsat}
		Let \(M\) be an \(R\)-module. Then the following statements hold.
		\begin{enumerate}
			\item If \(S \subseteq T\) are multiplicative subsets of \(R\) and \(M\) is a \(u\)-\(S\)-multiplication module, then \(M\) is a \(u\)-\(T\)-multiplication module.
			\item \(M\) is a \(u\)-\(S\)-multiplication module if and only if \(M\) is a \(u\)-\(S^{*}\)-multiplication module, where \(S^{*}\) is the saturation of \(S\).
		\end{enumerate}
	\end{proposition}
	
	\begin{proposition} \cite[Proposition 3]{Qi2025}
		\label{prop:usmultmodisoinv}
		Let \(M\) and \(M'\) be \(R\)-modules. Suppose \(M\) is \(u\)-\(S\)-isomorphic to \(M'\). Then \(M\) is a \(u\)-\(S\)-multiplication module if and only if \(M'\) is a \(u\)-\(S\)-multiplication module.
	\end{proposition}
	
	The property is also preserved under \(u\)-\(S\)-epimorphisms under suitable conditions.
	
	\begin{proposition} \cite[Proposition 4]{Qi2025}
		\label{prop:usmultmodepi}
		Let \(M\) and \(M'\) be \(R\)-modules. Suppose that \(S\) is a multiplicative subset of \(R\) and \(f:M \to M'\) is a \(u\)-\(S\)-epimorphism. If \(M\) is a \(u\)-\(S\)-multiplication module, then \(M'\) is a \(u\)-\(S\)-multiplication module. Conversely, suppose that \(M'\) is an \(S\)-multiplication module and \(t\Ker(f) = 0\) for some \(t \in S\); then \(M\) is a \(u\)-\(S\)-multiplication module.
	\end{proposition}
	
	Localization also preserves the property.
	
	\begin{proposition} \cite[Proposition 5]{Qi2025}
		\label{prop:usmultmodloc}
		Let \(R\) be a commutative ring and let \(S\) and \(T\) be multiplicative subsets of \(R\). Set \(\tilde{S} = \{\frac{s}{1} \in T^{-1}R \mid s \in S\}\). Suppose \(M\) is a \(u\)-\(S\)-multiplication \(R\)-module. Then \(T^{-1}M\) is a \(u\)-\(\tilde{S}\)-multiplication \(T^{-1}R\)-module.
	\end{proposition}
	
	When \(S\) satisfies the maximal multiple condition, the property is equivalent to being an \(S\)-multiplication module and to having a multiplication module after localization.
	
	\begin{proposition} \cite[Proposition 6]{Qi2025}
		\label{prop:usmultmodmmc}
		Let \(M\) be an \(R\)-module and let \(S\) be a multiplicative subset of \(R\) satisfying the maximal multiple condition. Then the following statements hold:
		\begin{enumerate}
			\item \(M\) is a \(u\)-\(S\)-multiplication module.
			\item \(M\) is an \(S\)-multiplication module.
			\item \(S^{-1}M\) is a multiplication \(S^{-1}R\)-module.
		\end{enumerate}
	\end{proposition}
	
	A connection with \(u\)-\(S\)-Noetherian modules is given: over a \(u\)-\(S\)-Noetherian ring, a \(u\)-\(S\)-multiplication module is \(u\)-\(S\)-Noetherian.
	
	\begin{proposition} \cite[Proposition 7]{Qi2025}
		\label{prop:usmultmodnoeth}
		Let \(R\) be a \(u\)-\(S\)-Noetherian ring and let \(M\) be a \(u\)-\(S\)-multiplication \(R\)-module. Then \(M\) is a \(u\)-\(S\)-Noetherian \(R\)-module.
	\end{proposition}
	
	The behavior under idealization is also characterized.
	
	\begin{theorem} \cite[Theorem 2]{Qi2025}
		\label{thm:usmultmodidealization}
		Let \(M\) be an \(R\)-module, let \(N\) be a submodule of \(M\), and let \(S\) be a multiplicative subset of \(R\). Then the following statements are equivalent.
		\begin{enumerate}
			\item \(N\) is a \(u\)-\(S\)-multiplication \(R\)-module.
			\item \(0(+)N\) is a \(u\)-\(S(+0)\)-multiplication ideal of \(R(+)M\).
			\item \(0(+)N\) is a \(u\)-\(S(+M)\)-multiplication ideal of \(R(+)M\).
		\end{enumerate}
	\end{theorem}
	
	A ring is a \(u\)-\(S\)-multiplication ring if it is a \(u\)-\(S\)-multiplication module over itself.
	
	\begin{definition} \cite[Definition 3]{Qi2025}
		\label{def:usmultring}
		Let \(R\) be a ring and let \(S\) be a multiplicative subset of \(R\). Then \(R\) is called a \textbf{\(u\)-\(S\)-multiplication} ring (with respect to \(s\)) if there exists \(s \in S\) such that each ideal of \(R\) is a \(u\)-\(S\)-multiplication with respect to \(s\), equivalently, if there exists \(s \in S\) such that, for each pair of ideals \(J \subseteq K\) of \(R\), there exists an ideal \(I\) of \(R\) satisfying \(sJ \subseteq IK \subseteq J\).
	\end{definition}
	
	Every multiplication ring is \(u\)-\(S\)-multiplication for any \(S\).
	
	\begin{corollary}\cite[Corollary 1]{Qi2025}
		\label{cor:multringisusmult}
		Every multiplication ring is a \(u\)-\(S\)-multiplication ring.
	\end{corollary}
	
	The property for a finite direct product decomposes into the properties for the factors.
	
	\begin{proposition} \cite[Proposition 9]{Qi2025}
		\label{prop:usmultringproduct}
		Let \(R = R_1 \times R_2\) and \(S = S_1 \times S_2\). Then \(R\) is a \(u\)-\(S\)-multiplication ring if and only if \(R_1\) is a \(u\)-\(S_1\)-multiplication ring and \(R_2\) is a \(u\)-\(S_2\)-multiplication ring.
	\end{proposition}
	
	When \(S\) satisfies the maximal multiple condition, the \(S\)-multiplication and \(u\)-\(S\)-multiplication properties coincide.
	
	\begin{proposition} \cite[Proposition 10]{Qi2025}
		\label{prop:usmultringmmc}
		Let \(S\) be a multiplicative subset of \(R\) that satisfies the maximal multiple condition. Then \(R\) is an \(S\)-multiplication ring if and only if \(R\) is a \(u\)-\(S\)-multiplication ring.
	\end{proposition}
	
	Localizing a \(u\)-\(S\)-multiplication ring at the uniform element yields a classical multiplication ring.
	
	\begin{proposition} \cite[Proposition 11]{Qi2025}
		\label{prop:usmultringloc}
		Suppose \(R\) is a \(u\)-\(S\)-multiplication ring. Then there is an \(s \in S\) such that \(R_s\) is a multiplication ring.
	\end{proposition}
	
	An example shows that a ring where every ideal is an \(S\)-multiplication ideal need not be a \(u\)-\(S\)-multiplication ring, and that the converse of Proposition \ref{prop:usmultringloc} is false.
	
	\begin{example}\cite[Example 3]{Qi2025}
		\label{ex:sidealusnotusmultring}
		Let \(D\) be an integral domain such that \(D_s\) is not a Dedekind domain for any \(0 \neq s \in D\) (e.g., \(D = k[x_1, x_2, \dots]\)). Set \(S = D - \{0\}\). Then every ideal of \(D\) is a \(u\)-\(S\)-multiplication ideal (since for any nonzero \(K\), pick \(s \in K\), then \(sK \subseteq sR \subseteq K\)), so \(D\) is an \(S\)-multiplication ring. However, \(D\) is not a \(u\)-\(S\)-multiplication ring by Proposition \ref{prop:usmultringloc}.
	\end{example}
	
	A remark clarifies that a ring where \(R_s\) is a multiplication ring for some \(s\) need not be \(u\)-\(S\)-multiplication.
	
	\begin{remark} \cite[Remark 1]{Qi2025}
		\label{rem:multlocnotusmult}
		Let \(D\) be a valuation domain with valuation group \(\mathbb{Z} \times \mathbb{Z}\) and maximal ideal generated by \(s\). Set \(S = D \setminus \{0\}\). Then \(D_s\) is a discrete valuation domain, hence a multiplication ring, but \(D\) is not \(u\)-\(S\)-multiplication by the argument in Example \ref{ex:sidealusnotusmultring}.
	\end{remark}
	
	A local characterization of classical multiplication rings is given.
	
	\begin{theorem}\cite[Theorem 3]{Qi2025}
		\label{thm:multringlocalus}
		Let \(R\) be a ring. Then the following statements are equivalent:
		\begin{enumerate}
			\item \(R\) is a multiplication ring.
			\item \(R\) is a \(u\)-\(\mathfrak{p}\)-multiplication ring for each \(\mathfrak{p} \in \Spec(R)\).
			\item \(R\) is a \(u\)-\(\mathfrak{m}\)-multiplication ring for each \(\mathfrak{m} \in \Max(R)\).
		\end{enumerate}
	\end{theorem}
	
	Finally, the property is studied for idealizations, showing that \(R(+)M\) being a \(u\)-\(S(+)M\)-multiplication ring implies that \(R\) is \(u\)-\(S\)-multiplication and every submodule of \(M\) is \(u\)-\(S\)-multiplication.
	
	\begin{proposition} \cite[Proposition 12]{Qi2025}
		\label{prop:usmultringidealization}
		Let \(R\) be a ring, let \(M\) be an \(R\)-module, and let \(S\) be a multiplicative subset of \(R\). Suppose \(R(+)M\) is a \(u\)-\(S(+)M\)-multiplication ring with respect to some \((s,m) \in S(+)M\). Then \(R\) is a \(u\)-\(S\)-multiplication ring with respect to \(s\), and each submodule of \(M\) is a \(u\)-\(S\)-multiplication \(R\)-module with respect to \(s\).
	\end{proposition}
	
	\subsection{Rings with Uniformly S-Noetherian Spectrum}
	
	The concept of a ring having a Noetherian spectrum is a topological condition on the set of prime ideals. The \(S\)-analogue was introduced by Hamed. Its uniform version, requiring a single element to control the radical finiteness, was studied by Guesmi and Hamed, and further by Zhang, Hamed, and Kim.
	
	\begin{definition} \cite[Definition 3.1]{Guesmi2025}
		\label{def:usnoethspec}
		Let \(R\) be a ring and \(S\) a multiplicative subset of \(R\). \(R\) is said to have a \textbf{uniformly \(S\)-Noetherian spectrum} (or \(u\)-\(S\)-Noetherian spectrum) if there exists \(s \in S\) such that every ideal of \(R\) is radically \(S\)-finite with respect to \(s\). That is, for every ideal \(I\) of \(R\), there exists a finitely generated subideal \(J \subseteq I\) such that \(sI \subseteq \sqrt{J} \subseteq \sqrt{I}\).
	\end{definition}

	The following result shows that when \(S\) satisfies the maximal multiple condition, the \(S\)-Noetherian spectrum and \(u\)-\(S\)-Noetherian spectrum properties coincide.
	
	\begin{proposition} \cite[ Proposition 2.1]{Zhang2025c}
		\label{prop:usnoethspecmmc}
		Let \(R\) be a commutative ring and \(S\) a multiplicative subset of \(R\) satisfying the maximal multiple condition. Then \(R\) is a ring with \(u\)-\(S\)-Noetherian spectrum if and only if \(R\) is a ring with \(S\)-Noetherian spectrum.
	\end{proposition}
	
	A key lemma is that localizing a ring with \(u\)-\(S\)-Noetherian spectrum at the uniform element yields a ring with Noetherian spectrum.
	
	\begin{lemma} \cite[Lemma 3.14]{Guesmi2025}
		\label{lem:usnoethspecloc}
		Let \(R\) be a ring with \(u\)-\(S\)-Noetherian spectrum with respect to some \(s \in S\). Then \(R_s\) is a ring with Noetherian spectrum.
	\end{lemma}
	
	An example shows that a ring can have an \(S\)-Noetherian spectrum without having a \(u\)-\(S\)-Noetherian spectrum.
	
	\begin{example} \cite[Example 3.15]{Guesmi2025}
		\label{ex:snoethspecnotusnoethspec}
		Let \(R = k[x_1, x_2, \dots]\) be a polynomial ring in countably many variables over a field \(k\), and set \(S = R \setminus \{0\}\). Then \(R\) is an \(S\)-Noetherian ring (and hence has an \(S\)-Noetherian spectrum), but it does not have the \(u\)-\(S\)-Noetherian spectrum property.
	\end{example}
	
	A characterization of rings with \(u\)-\(S\)-Noetherian spectrum, under the assumption that \(S\) consists of non-zero-divisors, involves the localization \(R_s\) and a compatibility condition.
	
	\begin{theorem}\cite[Theorem 2.4]{Zhang2025c}
		\label{thm:usnoethspecchar}
		Let \(R\) be a commutative ring and \(S\) a multiplicative subset consisting entirely of non-zero-divisors of \(R\). Then \(R\) has the \(u\)-\(S\)-Noetherian spectrum property if and only if there exists \(s \in S\) such that \(R_s\) has a Noetherian spectrum and
		\[
		\sqrt{I} R_s \cap R = \sqrt{I} :_R s
		\]
		for every finitely generated ideal \(I\) of \(R\).
	\end{theorem}
	
	In such rings, every radical ideal is a finite intersection of \(s\)-prime ideals.
	
	\begin{theorem} \cite[Theorem 3.1]{Zhang2025c}
		\label{thm:usnoethspecradicalint}
		Let \(R\) be a commutative ring and \(S\) a multiplicative subset of \(R\). If \(R\) has a \(u\)-\(S\)-Noetherian spectrum with respect to \(s\), then each radical ideal of \(R\) is a finite intersection of \(s\)-prime ideals.
	\end{theorem}
	
	Moreover, minimal \(s\)-prime ideals exist and every radical ideal is a finite intersection of them.
	
	\begin{lemma} \cite[Lemma 3.2]{Zhang2025c}
		\label{lem:minimalusprime}
		Let \(I\) be a proper ideal of \(R\), and let \(P\) be an \(s\)-prime ideal of \(R\) for some \(s \in S\) such that \(I \subseteq P\). Then there exists an \(s\)-prime ideal \(Q\) of \(R\) that is minimal among the \(s\)-prime ideals containing \(I\).
	\end{lemma}
	
	\begin{theorem} \cite[Theorem 3.3]{Zhang2025c}
		\label{thm:usnoethspecminimalint}
		Let \(R\) be a commutative ring, and let \(S\) be a multiplicative subset of \(R\). If \(R\) has a \(u\)-\(S\)-Noetherian spectrum with respect to \(s\), then each radical ideal of \(R\) is a finite intersection of minimal \(s\)-prime ideals.
	\end{theorem}
	
	When \(S\) consists of units, we recover the classical result.
	
	\begin{corollary} \cite[Corollary 3.4]{Zhang2025c}
		\label{cor:noethspecminimalint}
		If \(R\) has a Noetherian spectrum, then every radical ideal is a finite intersection of minimal prime ideals.
	\end{corollary}
	
	The property is preserved under various ring constructions.
	
	\begin{proposition} \cite[Proposition 4.1]{Zhang2025c}
		\label{prop:usnoethspecloc}
		Let \(R\) be a ring and \(S\) a multiplicative subset of \(R\). If \(R\) has \(u\)-\(S\)-Noetherian spectrum with respect to \(s \in S\), then \(U^{-1}R\) has \(u\)-\(U^{-1}S\)-Noetherian spectrum with respect to \(\frac{s}{1}\).
	\end{proposition}
	
	Flat overrings also inherit the property.
	
	\begin{theorem} \cite[Theorem 4.2]{Zhang2025c}
		\label{thm:usnoethspecflatoverring}
		Let \(R\) be a ring and \(S\) a multiplicative subset of \(R\). If \(R\) has \(u\)-\(S\)-Noetherian spectrum, then any flat overring of \(R\) also has \(u\)-\(S\)-Noetherian spectrum (identifying \(S\) with its image inside the flat overring of \(R\)).
	\end{theorem}
	
	The Hilbert basis theorem holds for this property: \(R[X]\) has \(u\)-\(S\)-Noetherian spectrum iff \(R\) does.
	
	\begin{theorem} \cite[Theorem 4.3]{Zhang2025c}
		\label{thm:usnoethspecHilbert}
		Let \(R\) be a commutative ring and \(S\) a multiplicative subset of \(R\). Then \(R\) has the \(u\)-\(S\)-Noetherian spectrum property if and only if the polynomial ring \(R[X]\) has the \(u\)-\(S\)-Noetherian spectrum property.
	\end{theorem}
	
	A remark notes that the power series ring may fail to have this property, even if \(R\) does.
	
	\begin{remark} \cite[Remark 4.4]{Zhang2025c}
		\label{rem:psnotnoethspec}
		By \cite{Ribenboim1985}, there exists a commutative ring \(R\) with Noetherian spectrum such that the power series ring \(R[[X]]\) does not have Noetherian spectrum. Consequently, even if \(R\) has \(u\)-\(S\)-Noetherian spectrum, \(R[[X]]\) may fail to have \(u\)-\(S\)-Noetherian spectrum.
	\end{remark}
	
	The property is also studied for fiber products (pullbacks) and amalgamations.
	
	\begin{theorem} \cite[Theorem 4.5]{Zhang2025c}
		\label{thm:usnoethspecpullback}
		Let \((T, \mathcal{M})\) be a (quasi-)local integral domain, \(K\) its residue field, \(\phi: T \to K\) the canonical surjection, \(D\) a subring of \(K\), and \(R = \phi^{-1}(D) \cong T \times_K D\). Let \(S_D\) (resp., \(S_T\)) be a multiplicative subset of \(D\) (resp., \(T\)), and \(S = \{s \in R \mid \lambda'(s) \in S_T \text{ and } \phi'(s) \in S_D\}\). Then \(R\) has \(u\)-\(S\)-Noetherian spectrum if and only if \(T\) has \(u\)-\(S_T\)-Noetherian spectrum and \(D\) has \(u\)-\(S_D\)-Noetherian spectrum.
	\end{theorem}
	
	Finally, the property is equivalent for a ring and several of its important ring extensions: the polynomial ring, Serre's conjecture ring, Nagata ring, and Anderson ring.
	
	\begin{theorem} \cite[Theorem 5.2]{Zhang2025c}
		\label{thm:usnoethspecvarious}
		Let \(R\) be a commutative ring with identity, and let \(S\) be a multiplicative subset of \(R\). The following statements are equivalent:
		\begin{enumerate}
			\item \(R\) has the \(u\)-\(S\)-Noetherian spectrum property.
			\item \(R[X]\) has the \(u\)-\(S\)-Noetherian spectrum property.
			\item \(R[X]_U\) (Serre's conjecture ring) has the \(u\)-\(S\)-Noetherian spectrum property.
			\item \(R[X]_N\) (Nagata ring) has the \(u\)-\(S\)-Noetherian spectrum property.
			\item \(R[X]_A\) (Anderson ring) has the \(u\)-\(S\)-Noetherian spectrum property.
		\end{enumerate}
	\end{theorem}
	
	When \(S\) consists of units, this recovers a classical result.
	
	\begin{corollary} \cite[Corollary 5.3]{Zhang2025c}
		\label{cor:noethspecvarious}
		The following statements are equivalent for a commutative ring \(R\):
		\begin{enumerate}
			\item \(R\) has Noetherian spectrum.
			\item \(R[X]\) has Noetherian spectrum.
			\item \(R[X]_U\) has Noetherian spectrum.
			\item \(R[X]_N\) has Noetherian spectrum.
			\item \(R[X]_A\) has Noetherian spectrum.
		\end{enumerate}
	\end{corollary}
	
	A sufficient condition for the semigroup ring \(R[\Gamma]\) to have the property is given.
	
	\begin{theorem} \cite[Theorem 5.4]{Zhang2025c}
		\label{thm:usnoethspecsemigroup}
		Let \(R\) be a commutative ring with identity, \(S\) a multiplicative subset of \(R\), and \(\Gamma\) a monoid. If \(R\) has the \(u\)-\(S\)-Noetherian spectrum property and \(\Gamma\) is finitely generated, then \(R[\Gamma]\) also has the \(u\)-\(S\)-Noetherian spectrum property.
	\end{theorem}
	
	For a numerical semigroup, the property is equivalent for the ring and its semigroup ring.
	
	\begin{corollary} \cite[Corollary 5.5]{Zhang2025c}
		\label{cor:usnoethspecnumsg}
		Let \(R\) be a commutative ring with identity, \(S\) a subset of \(R\), and \(\Gamma\) a proper numerical semigroup. Then \(R\) has the \(u\)-\(S\)-Noetherian spectrum property if and only if \(R[\Gamma]\) has the \(u\)-\(S\)-Noetherian spectrum property.
	\end{corollary}

\end{document}